\documentclass[11pt]{amsart}
\usepackage{fullpage,amssymb,mathtools,hyperref,pifont,microtype,mathrsfs,tikz,comment,verbatim}
\usepackage[T1]{fontenc}

\definecolor{blue}{rgb}{0,0,1}
\definecolor{red}{rgb}{1,0,0}
\definecolor{green}{rgb}{0,.6,.2}
\definecolor{purple}{rgb}{1,0,1}
\definecolor{teal}{rgb}{0,0.7,0.8}

\long\def\red#1\endred{\textcolor{red}{#1}}
\long\def\blue#1\endblue{\textcolor{blue}{#1}}
\long\def\purple#1\endpurple{\textcolor{purple}{ #1}}
\long\def\green#1\endgreen{\textcolor{green}{#1}}
\long\def\teal#1\endteal{\textcolor{teal}{#1}}

\newcommand{\xmark}{{\color{red}\ding{55}}}
\newcommand{\C}{\mathbb{C}}
\newcommand{\R}{\mathbb{R}}
\newcommand{\Q}{\mathbb{Q}}
\newcommand{\Z}{\mathbb{Z}}
\newcommand{\N}{\mathbb{N}}
\renewcommand{\H}{\mathbb{H}}
\newcommand{\E}{\mathcal{E}}

\newcommand{\calH}{\mathcal{H}}
\newcommand{\calP}{\mathcal{P}}
\newcommand{\cL}{\mathcal{L}}
\newcommand{\calR}{\mathcal{R}}
\newcommand{\K}{\mathcal{K}}
\newcommand{\W}{\mathcal{W}}
\newcommand{\F}{\mathbb{F}}

\newcommand{\scrE}{\mathscr{E}}

\newcommand{\cI}{\mathcal{I}}
\newcommand{\mfm}{\mathfrak{m}}
\newcommand{\vecepsilon}{{\pmb{\epsilon}}}
\newcommand{\veclambda}{{\pmb{\lambda}}}

\DeclareMathOperator{\Tr}{Tr}
\DeclareMathOperator{\Res}{Res}
\DeclareMathOperator{\End}{End}
\DeclareMathOperator{\Hom}{Hom}
\DeclareMathOperator{\dist}{dist}
\DeclareMathOperator{\lcm}{lcm}
\DeclareMathOperator{\ord}{ord}

\DeclareMathOperator{\rad}{rad}

\newtheorem{theorem}{Theorem}
\newtheorem{lemma}[theorem]{Lemma}
\newtheorem{proposition}[theorem]{Proposition}
\newtheorem{corollary}[theorem]{Corollary}
\theoremstyle{remark}

\numberwithin{theorem}{section}
\numberwithin{equation}{section}

\title{Towards a twistless converse theorem}
\author{Andrew R. Booker}
\author{Min Lee}
\address{School of Mathematics, University of Bristol, Bristol, United Kingdom}

\author{Nina Zubrilina}
\address{Department of Mathematics, Courant Institute of Mathematical Sciences, New York University, 251 Mercer Street, New York, NY 10012-1185, United States}

\begin{document}

\begin{abstract}
We prove a finiteness theorem for Dirichlet series satisfying the Euler product, analytic continuation, and functional equation expected of a holomorphic newform. For any fixed level $N$, even weight $k$, and exponent $\sigma$, only finitely many such series have Dirichlet coefficients $\ll n^\sigma$. No hypotheses on twists by Dirichlet characters are assumed. The proof combines a combinatorial spanning theorem for Ramanujan sums with an operator-theoretic argument that uses two distinct primes to establish discreteness of the possible Hecke eigenvalues.
\end{abstract}
\maketitle
\tableofcontents

\section{Introduction}

For $N,k\in\N$ with $k$ even, let $\E(N,k)$ be the set of complex-valued sequences $(f_n)_{n\in\N}$ satisfying the following hypotheses:
\begin{itemize}
\item the $L$-series of $f$ converges absolutely in some right half-plane and has a degree $2$ Euler product:
\begin{equation}\label{e:EP}
L_f(s) \coloneq \sum_{n=1}^\infty f_nn^{-s}
= \prod_p\frac1{1-f_p p^{-s} + \chi_0(p)p^{-2s}}
\quad\text{for }\Re(s)\gg_f1,
\end{equation}
where $\chi_0$ is the trivial character modulo $N$, and for each prime $p\mid N$,
\begin{equation}\label{e:AL}
f_p^2=\begin{cases}
0 & \text{ if } p^2\mid N, \\ p^{-1} & \text{ if } p\parallel N;
\end{cases}
\end{equation}
\item the completed $L$-function 
\[\Lambda_f(s) \coloneq N^{\frac12(s-\frac{1}{2})} \Gamma_\C\!\left(s+\frac{k-1}{2}\right) L_f(s),
\quad\Gamma_\C(s)\coloneq2(2\pi)^{-s}\Gamma(s),\]
extends to an entire function of finite order and satisfies the functional equation
\begin{equation}\label{e:FE}
\Lambda_f(s) = \epsilon_f \Lambda_f(1-s) \qquad\text{ for some } \epsilon_f \in \{\pm1\}.
\end{equation}
\end{itemize}

The above hypotheses hold when $f_n$ are the normalized Fourier coefficients of a Hecke eigenform in $M_k^\circ(\Gamma_0(N))$, the space of forms in $M_k^{\rm new}(\Gamma_0(N))$ that vanish at the cusps $0$ and $\infty$.
Farmer, Koutsoliotas, and Lemurell conjectured \cite[\S5, (1)]{FKL} that the converse is true: given $(f_n)_{n\in \N} \in\E(N,k)$, the function $f(z) \coloneq \sum_{n=1}^\infty f_n n^{\frac{k-1}{2}} e^{2\pi inz}$ is a modular form of weight $k$ and level $N$. For sequences satisfying the Ramanujan hypothesis $f_n\ll n^{o(1)}$, this is a special case of the automorphy conjecture for the Selberg class \cite[Conjecture~6.2]{Perelli}.

The first converse theorem is due to Hecke \cite{Hecke}, who proved that for $N\le4$, the functional equation is sufficient to show that $(f_n)$ comes from a modular form. Hecke's theorem does not assume that the $L$-series factors into an Euler product. As Hecke observed, the proof does not generalize to larger $N$ without additional hypotheses: the vector space of sequences satisfying the functional equation \eqref{e:FE} is infinite-dimensional for $N \geq 5$. Weil \cite{Weil} remedied this by further assuming functional equations for twists of the $L$-function by Dirichlet characters of conductor coprime to $N$, obtaining modularity for all $N$. One inconvenient feature of this result is that the number of character twists required to conclude modularity grows with $N$ \cite{S19}. Another method for showing modularity without twists is that of Dimitrov \cite{Dimitrov}, who has proven a converse theorem that requires only one functional equation and assumes integrality of the unnormalized coefficients $f_n n^{\frac{k-1}{2}}$. Finally, the approach that we will pursue in this paper is via multiplicativity, as in \cite{FKL}. Assuming that the $L$-function satisfies an Euler product of the right shape, \cite{CF,CFOS,BCOZ,Zhang26} extended Hecke's ideas to prove converse theorems for $N\in\{5,\ldots,18\}\cup\{20,23,24,32,36\}$. 

One non-trivial consequence of the conjecture in \cite{FKL} is that the set $\E(N,k)$ is \emph{finite}. Our first result establishes this finiteness for sequences of a fixed polynomial growth rate. For $\sigma\in\R$, let
\[
\E_\sigma(N,k)\coloneq\bigl\{(f_n)_{n\in\N}\in\E(N,k):f_n=O(n^\sigma)\bigr\}.
\]
For instance, the conjecture implies that $\#\E_\sigma(N,k)=\dim S_k^{\rm new}(\Gamma_0(N))$ for any $\sigma\in(0,\frac{k-1}{2})$.
\begin{theorem}\label{thm:main1}
Let $N,k\in\N$ with $k$ even, and let $\sigma\in\R$. Then $\E_\sigma(N,k)$ is finite.
\end{theorem}
\noindent
For a large class of levels, we can remove the growth condition and show that the count is asymptotically the one predicted by the conjecture. Let $\rad(N)\coloneq\prod_{p\mid N}p$ denote the radical of $N$.
\begin{theorem}\label{thm:main-asymptotic}
Suppose $\rad(N)>\sqrt{N}$. Then $\E(N,k)$ is finite and
\[\#\E(N,k)=\dim M_k^\circ(\Gamma_0(N))+
O_N\!\left(\frac{k}{\log\log{2k}}\right).\]
\end{theorem}
\noindent
Finally, we obtain the full conjecture for some additional small levels:
\begin{theorem}\label{thm:small-levels}
Suppose $N\in\{1,\ldots,36\}\setminus\{29,31\}$. Then $f\in M_k^\circ(\Gamma_0(N))$ for every $(f_n)\in\E(N,k)$. In particular, $\#\E(N,k)=\dim M_k^\circ(\Gamma_0(N))$.
\end{theorem}
The bulk of the paper is devoted to the proof of Theorem~\ref{thm:main1}. Along the way, we prove the finiteness assertion of Theorem~\ref{thm:main-asymptotic} in the form $\#\E(N,k)=O_N(k)$ (Corollary~\ref{cor:main-dividingN} and the remark following it). We do not reproduce the proofs of the asymptotic formula in Theorem~\ref{thm:main-asymptotic} or Theorem~\ref{thm:small-levels} here. They are considerably longer and were written by large language models; complete proofs are available at \cite{github}.

The growth exponent in Theorem~\ref{thm:main1} and the condition in Theorem~\ref{thm:main-asymptotic} play the following roles. For $p\mid N$, the Euler product determines $f_p$ up to finitely many values, since $f_p^2\in\{0,p^{-1}\}$ by \eqref{e:AL}. When $\rad(N)>\sqrt N$, these primes suffice for the spanning argument (part~(2) of Theorem~\ref{thm:spanning}), so no Hecke eigenvalue is left undetermined and no growth hypothesis is needed. In general the argument needs a prime $p\nmid N$ with $p>\sqrt N$, for which $f_p$ is an unknown Hecke eigenvalue; its possible values form a discrete set (Proposition~\ref{prop:point-discrete}), and the growth exponent $\sigma$ serves only to confine $f_p$ to a bounded set.

\subsection*{Remarks}
\begin{itemize}
\item Here $M_k^{\rm new}(\Gamma_0(N))$ denotes the space of $f\in M_k(\Gamma_0(N))$ such that both $f$ and its Fricke image $f|W_N$ have trace $0$ from $\Gamma_0(N)$ to $\Gamma_0(N/p)$ for every prime $p\mid N$. The space $M_k^\circ(\Gamma_0(N))$ is spanned by the cuspidal subspace $S_k^{\rm new}(\Gamma_0(N))\subseteq M_k^{\rm new}(\Gamma_0(N))$ together with the Eisenstein series
\[
\sum_{d\mid r}\mu(d)\chi(d)d\sum_{n=1}^\infty\sum_{m\mid n}\chi(n/m)\overline{\chi}(m)m^{k-1}e^{2\pi idnz},
\]
where $N=q^2r$, $q>1$, $r$ is squarefree and coprime to $q$, $r=1$ unless $k=2$, and $\chi$ is a primitive Dirichlet character mod $q$ satisfying $\chi(p)\in\{\pm1\}$ for every prime $p\mid r$. In particular, $M_k^\circ(\Gamma_0(N))=S_k^{\rm new}(\Gamma_0(N))$ if $k>2$ and $N$ is not a square.
\item One key technical combinatorial ingredient of the proof of Theorem~\ref{thm:main1} is a spanning statement for certain explicit vectors with entries in $\{0,1\}$: see Theorem~\ref{thm:spanning} and its proof in \S\ref{s:spanning-proof}. This proof was conceived with the assistance of large language models.
\item Theorems~\ref{thm:main1}--\ref{thm:small-levels} have been formalized in Lean; formalizations are available at \cite{github}.
\item The proof of Theorem~\ref{thm:main1} can be adapted to omit the hypothesis on $f_p$ for $p\mid N$, and to allow non-trivial nebentypus characters and odd weights. We have stated the three main theorems for the common class $\E(N,k)$ for convenience.
\item A version of Theorem~\ref{thm:main-asymptotic} holds for Maass forms of weight $0$ and squarefree level: replacing $\Gamma_\C(s+\frac{k-1}{2})$ in the functional equation by $\Gamma_\R(s+\delta+ir)\Gamma_\R(s+\delta-ir)$ for fixed $\delta\in\{0,1\}$, $r\in i(0,\frac12)\cup\R_{\ge0}$, we can conclude that the number of solutions is finite for each fixed $r$, and the total number with $r\in[t,t+\varepsilon(t)]$ satisfies the expected asymptotic as $t\to\infty$, provided $\varepsilon(t)>0$ decays sufficiently slowly. In particular, we recover Weyl's law for the cuspidal spectrum without appealing to any trace formula.
While we focus on the holomorphic case in this paper, we include
a writeup of the Maass form case in \cite{github}.
\item For prime $N$, \cite{github} proves a version of Weil's converse theorem using twists of modulus $q<\sqrt{N}$. By \cite[Theorem~1.1]{S19}, this is sharp up to a constant factor.
\item The passage from finiteness to modularity in general is left as an exercise for the interested LLM.
\end{itemize}

\subsection*{Acknowledgements and LLM disclosure}
The idea of using Euler products to generalize the converse theorem was introduced by Conrey and Farmer \cite{CF} in 1995. 
While we have pursued the Ramanujan-sum approach since 2015, our work was greatly accelerated by recent advancements in AI. Starting in April 2026, we used Anthropic's Claude to conduct hundreds of numerical and thought experiments, eventually arriving at the key lemma, Theorem~\ref{thm:spanning}, as a conjecture. 
We continued to work with Claude and later OpenAI's ChatGPT, exploring various proof strategies, before finally obtaining the proof on September 18th. 
That proof relied on recent results of Byszewski, Konieczny, and Krawczyk \cite{BKK} on automatic sequences that we learned about from GPT-6 Astra.
The proof presented in \S\ref{s:spanning-proof} is a later simplification, in which this input is replaced by classical results of Cobham \cite{Cobham} and Furstenberg \cite{Furstenberg}.

\subsection*{Paper outline}
We first explain the heuristic behind the Ramanujan-sum approach. In \S\ref{sec:spanning} we state the combinatorial spanning theorem, Theorem~\ref{thm:spanning}, and in \S\ref{s:spanning-proof} we prove it. In \S\ref{sec:main-proof} we use reconstruction and operator theory to prove that $\E_\sigma(N,k)$ is finite, completing the proof of Theorem~\ref{thm:main1}.

\subsection*{Guiding heuristic: equations vs unknowns}
We begin by describing a heuristic for the Euler product approach. We borrow heavily from the forthcoming paper \cite{BS}, which applies the same heuristic in another context.

The idea is to view the functional equation \eqref{e:FE} as a linear system with unknowns $(f_n)$, truncated at $n \leq X$ for a careful choice of $X$. We control these unknowns with \eqref{e:FE} and its $t$-aspect shifts $L_f(s+it)$; to that end, we use the approximate functional equation to interpret \eqref{e:FE} as a linear condition involving the first $\sim \sqrt{C}$ coefficients $f_n$, where $C = C(N, k, t)$ is the analytic conductor.

Functional equations for $L_f(s + it)$ for nearby choices of $t$ carry similar information. As we move up the critical line $\{\frac12+it:t\in\R\}$, each interval of width $2\pi$ contains about $\log{C}$ zeros of $L$, and random matrix theory philosophy posits that each such zero packet corresponds to a distinct random matrix in a prescribed group. Imagining that functional equations in each packet produce one distinct relation, we expect to get $\sim\frac{T}{2\pi}$ independent linear constraints coming from $t\in[0,T]$. For large $t$, $C \sim N(\frac{t}{2\pi})^2$, and since we can go up to $C = X^2$, the number of equations we get is $\frac{T}{2\pi} = X/\sqrt{N}$. This is smaller than the number of unknowns, $X$, unless $N$ is very small, which is consistent with Hecke's result.

For the twist by a primitive Dirichlet character $\chi\bmod{q}$, the corresponding linear relation at height $t$ has length $\sqrt{C_\chi}\sim\frac{tq\sqrt{N}}{2\pi}$, and we can go up to $t\le T_q=\frac{2\pi X}{q\sqrt{N}}$. Hence, writing $\varphi^*(q)$ for the number of primitive characters mod $q$, each of the approximately $\frac12\varphi^*(q)$ conjugate pairs of primitive characters $\chi\bmod{q}$ contributes about $\frac{T_q}{2\pi}=\frac{X}{q\sqrt{N}}$ equations, for a total of $\sum_{q\le X/\sqrt{N}}\frac{\varphi^*(q) X}{2q\sqrt{N}} \asymp \frac{X^2}{N},$ making the system highly overdetermined for large $X$.

Finally, an Euler product implies that only the $f_p$ for $p$ prime are unknown. This saves a $\log$ factor in the variable count, and the system becomes overdetermined once $\log X \gg \sqrt{N}$. The problem, however, is that the resulting equations are polynomial, not linear, in the unknowns $f_p$.

One way around this complication is by \emph{linearizing the Euler product}.
An Euler product relation can be thought of as a recipe for an imprimitive character ``twist'' of an $L$-function, in the following sense. For a modulus $q$, $L$ has additive character twists by $e^{2\pi i an/q}$ with $(a, q) = 1$ and multiplicative character twists by $\chi(n)$ for $\chi$ primitive mod $q$. The two are related by Fourier analysis on $\Z/q\Z$ and $(\Z/q\Z)^\times$, with one complication: for prime $q$, there are $q-1$ Farey fractions $a/q$ but only $q-2$ primitive $\chi\bmod{q}$.

The twist by the trivial character mod $q$ does not have a nice functional equation, but there is a fix: instead we can twist by its mod $q$ Fourier transform, the Ramanujan sum
\[c_q(n) = \sum_{\substack{a\bmod{q}\\ \gcd(a, q)=1}} e^{2\pi in \frac{a}{q}} \qquad\text{ for } n\in \Z.\]
For comparison, a primitive character $\chi\bmod{q}$ satisfies
\[\chi(n) = \frac{1}{\tau(\overline{\chi})}\sum_{a\bmod{q}} \overline{\chi}(a) e^{2\pi in \frac{a}{q}} \qquad\text{ for } n\in \Z,\]
where $\tau(\overline{\chi})\coloneq\sum_{a\bmod{q}}\overline{\chi}(a)e^{2\pi i\frac{a}{q}}$ is the Gauss sum; replacing $\overline{\chi}$ by the trivial character mod $q$, the sum on the right becomes $c_q(n)$.
As we show in Lemma~\ref{lem:Euler-Ramanujan}, when $L_f(s)$ has an Euler product as in \eqref{e:EP}, its twist by $c_q$ satisfies a functional equation of conductor $\ell_N(q)\coloneq\lcm(N,q^2)$,  similarly to a twist by a primitive character mod~$q$.

For each $q$, we can combine the $c_q$-twist with the $t$-aspect heuristic, and this gives a linear relation of length approximately $\frac{t\sqrt{\ell_N(q)}}{2\pi}$. Considering all $q$ and counting equations as above that involve unknowns $f_n$ for $n\le X$, we expect about
\[\sum_{\substack{q\in\N\\\ell_N(q)\le X^2}}
\frac{X}{\sqrt{\ell_N(q)}}
\ge\sum_{\substack{q\in\N\\Nq^2\le X^2}}\frac{X}{\sqrt{Nq^2}}
\ge\frac{X}{\sqrt{N}}\log\frac{X}{\sqrt{N}}
\]
equations. This again leads to the conclusion that the system is overdetermined by a log factor for large $X$, now in a linear system. Based on this, Bettin et al.\ \cite[Conjecture~1.1]{src} conjectured that the functional equations satisfied by Ramanujan-sum twists should suffice to imply modularity, and proved it for $N \in \{1, \ldots, 9\} \cup \{11, 15, 17, 23\}$.

\subsection{The spanning theorem}\label{sec:spanning}

Theorem~\ref{thm:main1} relies on a key combinatorial result about explicit vectors of length $Q$ with entries in $\{0,1\}$ spanning $Q$-dimensional space. 
We state the result in this subsection; its proof is given in \S\ref{s:spanning-proof}.

Fix a positive integer $N$ and recall $\ell_N(q) = \lcm(N, q^2)$ for $q\in \N$. 
By $\Q^{\Z/Q\Z}$ for a positive integer $Q$, we mean functions from $\Z/Q\Z$ to $\Q$, or equivalently, vectors of length $Q$ with entries in $\Q$, where the coordinates are considered modulo $Q$ for convenience. 
For $x \in \Q^{\Z/Q\Z}$, we denote the $n$-th coordinate of $x$ by $x(n)$. 
For $j\bmod{Q}$, we let $e_j$ denote the standard coordinate vector with entry $1$ at position $j$ and $0$ elsewhere, so $x = \sum_j x(j) e_j$. 
For $u \in \Z$ and $q \mid Q$, we define $R^Q_{q, u} \in \Q^{\Z/Q\Z}$ via
\begin{equation}\label{e:Rqu}
R^Q_{q, u} \coloneq \sum_{\substack{a\bmod{q} \\ \gcd(a, q)=1}} e_{u+a\frac{Q}{q}}
\end{equation}
and let
\begin{equation}\label{e:scrR_NQ_def}
\calR_N(Q) \coloneq \left\{R^Q_{q, u} : q\mid Q,\;u \in \Z\cap\left(-\frac{Q}{\sqrt{\ell_N(q)}}, \frac{Q}{\sqrt{\ell_N(q)}}-1\right)\right\}.
\end{equation}
We call the elements of $\calR_N(Q)$ \textit{legal rows}.
Equivalently, $\max\{|u|,|u+1|\}<Q/\sqrt{\ell_N(q)}$: the condition holds at both endpoints of the interval $[u/Q,(u+1)/Q]$, which is how it is used in \S\ref{sec:main-proof} (see \eqref{e:buffer}).
We also define the complete rows: for $q\mid Q$ and $u\in \Z$, 
\begin{equation}
C_{q, u}^Q \coloneq \sum_{a\bmod{q}} e_{u+a\frac{Q}{q}} = \sum_{d\mid q} R_{d, u}^Q. 
\end{equation}
We say that a finite set of primes $\calP$ is \emph{vanishing} for $N$ if there exists $Q\in(\prod_{p\in\calP}p)\N$ such that 
the vectors in $\calR_{N}(Q) \cup \bigl\{C_{p, u}^Q:\, p\in \calP, \, u\bmod{Q/p}\bigr\}$ span $\Q^{\Z/Q\Z}$ (equivalently, no nonzero $x\in\Q^{\Z/Q\Z}$ is orthogonal to all of them).

\begin{theorem}[Spanning]\label{thm:spanning}
Let $N$ be a positive integer. 
A finite set of primes $\calP$ is vanishing for $N$ if and only if for every $s\in\N$,
\begin{equation}\label{e:vanishing_P}
s>2 \text{ and } p\mid s \text{ for all }p\in\calP \implies s^2> \gcd(N,s^\infty)=\prod_{p\mid s} p^{\ord_p(N)}.
\end{equation}
In particular,
\begin{enumerate}
\item the empty set is vanishing for $N$ if and only if $N$ is divisible neither by $2^4$ nor by $p^2$ for any odd prime $p$;
\item the set of prime divisors of $N$ is vanishing for $N$ if and only if $N\mid8$ or $\rad(N) > \sqrt{N}$;
\item for every prime $p>\sqrt{N}$, the singleton set $\{p\}$ is vanishing for $N$.
\end{enumerate}
\end{theorem}

The condition in part~(1) holds exactly when $\Gamma_0(N)$ has one cusp of each denominator $d\mid N$: the cusps of $\Gamma_0(N)$ with denominator $d$ number $\varphi(\gcd(d,N/d))$, which equals $1$ exactly when $\gcd(d,N/d)\le2$.

In the table below, we record the minimum value of $Q$ for small values of $N$ when the empty set is vanishing for $N$; an entry \xmark{} means that it is not, by part~(1) of Theorem~\ref{thm:spanning}. 
See \cite{github} for details of the algorithms underlying these computations.

\begin{center}\scriptsize
\begin{tabular}{c|cccccccccc}
$N$ & 1 & 2 & 3 & 4 & 5 & 6 & 7 & 8 & 9 & 10 \\
$Q$ & 2 & 2 & 2 & 4 & 6 & 4 & 6 & 4 & \xmark & 10 \\\hline
$N$ & 11 & 12 & 13 & 14 & 15 & 16 & 17 & 18 & 19 & 20 \\
$Q$ & 44 & 4 & 104 & 12 & 20 & \xmark & 408 & \xmark & 570 & 20 \\\hline
$N$ & 21 & 22 & 23 & 24 & 25 & 26 & 27 & 28 & 29 & 30 \\
$Q$ & 84 & 132 & 1380 & 12 & \xmark & 312 & \xmark & 42 & 4872 & 30 \\\hline
$N$ & 31 & 32 & 33 & 34 & 35 & 36 & 37 & 38 & 39 & 40 \\
$Q$ & 10416 & \xmark & 264 & 1020 & 420 & \xmark & 15540 & 1596 & 780 & 90 \\\hline
$N$ & 41 & 42 & 43 & 44 & 45 & 46 & 47 & 48 & 49 & 50 \\
$Q$ & 34440 & 168 & 36120 & 990 & \xmark & 9660 & 78960 & \xmark & \xmark & \xmark \\\hline
$N$ & 51 & 52 & 53 & 54 & 55 & 56 & 57 & 58 & 59 & 60 \\
$Q$ & 4080 & 2730 & 133560 & \xmark & 7920 & 630 & 15960 & 36540 & 297360 & 210 \\\hline
$N$ & 61 & 62 & 63 & 64 & 65 & 66 & 67 & 68 & 69 & 70 \\
$Q$ & 614880 & 39060 & \xmark & \xmark & 32760 & 3300 & 1857240 & 10710 & 38640 & 1260 \\\hline
$N$ & 71 & 72 & 73 & 74 & 75 & 76 & 77 & 78 & 79 & 80 \\
$Q$ & 3936240 & \xmark & 8094240 & 512820 & \xmark & 131670 & 166320 & 5460 & 43797600 & \xmark \\\hline
$N$ & 81 & 82 & 83 & 84 & 85 & 86 & 87 & 88 & 89 & 90 \\
$Q$ & \xmark & 1704780 & 59819760 & 11550 & 942480 & 7747740 & 3483480 & 90090 & 384864480 & \xmark \\\hline
$N$ & 91 & 92 & 93 & 94 & 95 & & & & & \\
$Q$ & 3603600 & 2072070 & 7447440 & 25405380 & 13693680 & & & & &
\end{tabular}
\end{center}

Note that in order for $\calR_N(Q)$ to have at least $Q$ elements, $Q$ must be highly divisible, and hence large: 
the minimal spanning $Q$ has to satisfy $Q\ge\exp\exp\Bigl(N^{\frac12-\frac{2+o(1)}{\log\log{N}}}\Bigr)$ as $N\to\infty$. 
Our proof does not give a useful quantitative bound for $Q$.

Combined with the $c_q$-twisted functional equation, Theorem~\ref{thm:spanning} can be used to express $f(z)$ for $z$ near the real axis as a linear combination of values of $f$ at other arguments, and the restriction to legal rows guarantees that those arguments will have a larger imaginary part than $z$. 
Indeed, the $c_q$-twisted functional equation relates the sum $\sum_{a\in(\Z/q\Z)^\times}f(z+\frac aq)$ to the same sum evaluated at $-1/(\ell_N(q)z)$, and for $z=x+iy$ we have $\Im\bigl(-1/(\ell_N(q)z)\bigr)=y/(\ell_N(q)|z|^2)$, which exceeds $y$ when $\ell_N(q)|z|^2<1$. For small $y$ this holds when $|x|<1/\sqrt{\ell_N(q)}$; the row $R^Q_{q,u}$ is used at points with real part in $[u/Q,(u+1)/Q]$, and it is legal exactly when this interval lies in that range. 
A complete row $C^Q_{p,u}$ with $p\in\calP$ corresponds instead to the Euler factor at $p$: for $(f_n)\in\E(N,k)$,
\[\sum_{b\bmod p}f\Bigl(z+\frac bp\Bigr)=p^{\frac{k+1}{2}}f_p\,f(pz)-\chi_0(p)p^kf(p^2z),\]
and the arguments on the right have heights $p\Im z$ and $p^2\Im z$. 
The resulting reconstruction formula bounds each fixed Hecke eigenspace by finitely many Fourier coefficients. A compactness estimate and a comparison between two Hecke operators then show that only finitely many eigenvalues occur under a fixed growth bound. 

\section{Proof of Theorem \ref{thm:spanning}}\label{s:spanning-proof}

In this section we prove Theorem~\ref{thm:spanning}. 
The integer $N$ is fixed throughout the section. 
We write $\calP$ for the set in Theorem~\ref{thm:spanning}, reserving $\cL$ for a variable set of primes that may contain $\calP$. 
Recall that $\ell_N(q)=\lcm(N,q^2)$ and that the set $\calR_N(Q)$ of legal rows is defined as in \eqref{e:scrR_NQ_def}.
Note that if $d\mid q$ then $\ell_N(d) \mid \ell_N(q)$, so $R_{d, u}^Q$ is legal for any divisor $d$ of $q$ if $R_{q, u}^Q$ is legal. 
So we say the complete row $C_{q, u}^Q$ is legal if $R_{q, u}^Q$ is legal. 

Let $\F$ be a field. The rows $R^Q_{q,u}$ have entries in $\{0,1\}$, so they act on $\F^{\Z/Q\Z}$ by the pairing $\langle R,x\rangle=\sum_{n\bmod{Q}} R(n)x(n)$. 
Write
\begin{equation}\label{e:KQ_def}
\K(Q;\F) = \bigl\{x\in\F^{\Z/Q\Z}:\ \langle R,x\rangle=0\text{ for all }R\in\calR_N(Q)\bigr\}
\end{equation} 
for the common kernel of the legal rows. 

If $\ell\mid Q$, the \emph{trace} from $Q$ down to $Q/\ell$ is the map
\begin{equation}\label{sp:eq:trace}
\Tr_\ell^Q:\F^{\Z/Q\Z}\longrightarrow\F^{\Z/(Q/\ell)\Z},
\qquad (\Tr_\ell^Qx)(n) = 
\left<C_{\ell, n}^Q,x\right> = \sum_{j=0}^{\ell-1}x(n+jQ/\ell).
\end{equation}
Let $\cL$ be a finite set of primes. 
For every $Q$ divisible by all elements of $\cL$, put
\begin{equation}\label{sp:eq:augmented}
\W_{\cL}(Q;\F)=\K(Q;\F)\cap\bigcap_{\ell\in\cL}\ker\Tr_\ell^Q,
\qquad \W_{\emptyset}(Q;\F)=\K(Q;\F).
\end{equation} 
Since the rows $C_{\ell, u}^Q$ for $\ell\in \cL$ are not necessarily legal, adding these relations can make the space $\W_{\cL}(Q;\F)$ smaller. 
Recall that $\cL$ is vanishing if $\W_{\cL}(Q; \Q)=0$ for some $Q$ divisible by all elements of $\cL$. 

For sufficiency in Theorem \ref{thm:spanning}, we prove that $\W_{\calP}(Q;\Q)= 0$ for some $Q$, where $\calP$ satisfies \eqref{e:vanishing_P}.
We first find a vanishing set $\cL$ which contains $\calP$ and then remove primes from $\cL\setminus \calP$, one at a time. 
The proof rests on the following two propositions.

\begin{proposition}\label{sp:prop:initial}
Suppose that $N\ge5$. 
Let $\cL$ be a finite set of primes $p$ with $\gcd(p, 2N) = 1$ such that
\begin{equation}\label{sp:eq:rho}
\rho(\cL)\coloneq\prod_{p\in\cL}\bigl(1-p^{-1}\bigr)<\frac2{\sqrt N},
\end{equation}
and let $D=\prod_{p\in\cL}p$. Then $\W_{\cL}(D;\Q)= 0$, so $\cL$ is vanishing.
\end{proposition}

\begin{proposition}\label{sp:prop:two}
Suppose that $\calP$ satisfies \eqref{e:vanishing_P}, and let $\cL$ be a finite set of primes containing $\calP$. 
Let $c_1\ne c_2$ be primes not in $\cL$.
If $\cL\cup\{c_1\}$ and $\cL\cup\{c_2\}$ are vanishing, then $\cL$ is vanishing.
\end{proposition}

\begin{proof}[Proof of Theorem~\ref{thm:spanning}]
We first prove necessity.  
Let $N_s=\prod_{p\mid s}p^{\ord_p(N)}$.
Suppose $s>2$, every prime in $\calP$ divides $s$, and $s^2\le N_s$.
We claim that $\calP$ is not vanishing. 
At any length $Q$ divisible by $s$, choose a primitive $s$th root
of unity $\zeta$ and define
\[x(n) = \begin{cases}
\zeta^{\frac{n}{Q/s}} - \zeta^{-\frac{n}{Q/s}} & \text{ if } n\equiv 0\bmod{Q/s}, \\
0 & \text{ otherwise } 
\end{cases}
\qquad \text{ for } n\bmod{Q}.\]
This vector is nonzero and $x(0)=0$. 
For a complete row $C_{q, u}^{Q}$, we have $\left<C_{q, u}^Q, x\right>=0$ unless $\frac{Q}{\lcm(q, s)}\mid u$. 
When $\frac{Q}{\lcm(q, s)}\mid u$, we have 
\[\left<C_{q, u}^Q, x\right> = \sum_{a\bmod{q}} x(u+a Q/q) 
= \sum_{\substack{a\bmod{q}\\ \frac{s}{\gcd(q, s)} a \equiv -\frac{u}{Q/\lcm(q, s)} \bmod{\frac{q}{\gcd(q, s)}}}} \left(\zeta^{\frac{u+aQ/q}{Q/s}} - \zeta^{-\frac{u+aQ/q}{Q/s}}\right).\]

When $\gcd(q,s)=1$, the sum has a single term, which vanishes when $u=0$. If $u\neq0$, the divisibility condition above forces the row to be illegal. Indeed, $q\notin\calP$, and
since $u$ is divisible by $\frac{Q}{qs}$ and  
\[\lcm(q, s)^2 = q^2 s^2 \leq q^2N_s \leq \ell_N(q),\]
we have $|u| \geq \frac{Q}{\sqrt{\ell_N(q)}}$ if $u\neq 0$. 

When $\gcd(q, s) >1$, since $\sum_{b\bmod{\gcd(q, s)}} \zeta^{b\frac{s}{\gcd(q, s)}} =0$, we have $\left<C_{q, u}^Q, x\right>=0$.

Thus, $x\in \W_{\calP}(Q; \C)$ and $x\neq 0$.
If $\calP$ is vanishing, there exists $Q_1\in (\prod_{p\in \calP} p)\N$ such that $\W_{\calP}(Q_1, \Q)=0$. 
By Corollary \ref{cor:kernel_growingQ}, $\W_{\calP}(sQ_1, \Q)=0$ and also $\W_{\calP}(sQ_1, \C)=0$, which is a contradiction. 

For sufficiency, we assume \eqref{e:vanishing_P} for $\calP$. 
For $N\le4$, the $q=1$ rows at $Q=2$ for $N\in \{1,2,3\}$ and $Q=4$ for $N=4$, 
together in the last case with $R^4_{2,-1}$ and $R^4_{2,0}$, span directly.
If $\cL$ is vanishing, then so is $\cL\cup\calP$: by Corollary~\ref{cor:kernel_growingQ} we may take $Q$ divisible by every prime of $\calP$, and $\W_{\cL\cup\calP}(Q;\Q)\subseteq\W_{\cL}(Q;\Q)$.
For $N\ge5$, enumerate the primes outside $2N\prod_{p\in\calP}p$ as $p_1<p_2<\cdots$. 
The divergence of $\sum_p1/p$ permits a
choice $n\ge2$ such that $\prod_{i=2}^n(1-p_i^{-1})<2/\sqrt N$.
Every subset of $\{p_1,\ldots,p_n\}$ with at least $n-1$ elements is vanishing by Proposition~\ref{sp:prop:initial}; so is its union with $\calP$. 
Descending induction using
Proposition~\ref{sp:prop:two} now removes all the auxiliary primes.
Indeed, a subset of size at most $n-2$ has two distinct omitted
primes, and adjoining either gives a larger subset already handled.

Therefore \eqref{e:vanishing_P} is equivalent to $\calP$ being vanishing. 

For (1), if $16\mid N$ take $s=4$, and if an odd $p^2\mid N$ take $s=p$. 
Conversely, when neither divisibility holds, the factors of $s^2/N_s$ at odd primes dividing $s$ are at least $p$, and its factor at $2$ is at least $1/2$. 
Thus $s>2$ containing an odd prime has $s^2/N_s\ge3/2$; 
if $s$ is a power of $2$ then $s\ge4$ and $s^2/N_s\ge2$.
For (2), put $r=\rad(N)$. 
If $r>2$, the minimum eligible $s$ is $r$ and $N_s=N$, so the condition is exactly $r^2>N$. 
If $r=2$, the minimum is $s=4$ and the condition is $N<16$, namely $N\mid8$;
$N=1$ also satisfies the condition. 
(3) follows from
$s\ge p>\sqrt N\ge\sqrt{N_s}$.
\end{proof}

In the rest of the section, we prove Propositions \ref{sp:prop:initial} and \ref{sp:prop:two}.
The linear maps used throughout are introduced in \S\ref{sp:sec:algebra}, and Proposition \ref{sp:prop:initial} is proved in \S\ref{sp:sec:removal}; that proof uses only the trace conditions and the legal rows $R^Q_{1,u}=e_u$. 
The proof of Proposition~\ref{sp:prop:two} is laid out in \S\ref{sp:sec:stationary}--\S\ref{sp:sec:closure} and has four steps.
Conceptually, it works by climbing the two-parameter tower of lengths $Q=Q_0c_1^ac_2^b$, $a,b\geq0$: the traces go down the tower, the restrictions to residue classes compare an element at one level with elements at lower levels, and we study how the spaces $\W_{\cL}(Q;\F)$ and the coordinates of their elements behave as $a$ and $b$ grow.
\begin{enumerate}
\item \emph{Dimension stabilization} (\S\ref{sp:sec:stationary}). The hypothesis that $\cL\cup\{c_1\}$ and $\cL\cup\{c_2\}$ are vanishing is used as follows. Over a suitable finite field $\F$, and for every length $Q$ divisible by a suitable integer $B$, it says that for $c\in\{c_1,c_2\}$ no nonzero element of $\W_{\cL}(cQ;\F)$ has trace zero; that is, the trace $\Tr_c^{cQ}:\W_{\cL}(cQ;\F)\to\W_{\cL}(Q;\F)$ is injective (Lemma~\ref{lem:trace_injective}). So the dimension of $\W_{\cL}(Q;\F)$ cannot increase when $Q$ is multiplied by $c_1$ or $c_2$, and hence it is constant from some length $Q_0$ and up: for every $Q=Q_0c_1^ac_2^b$ with $a,b\geq0$, these traces are isomorphisms (Lemma~\ref{lem:stabilization}). This identifies all the spaces $\W_{\cL}(Q_0c_1^ac_2^b;\F)$ with the single finite space $V=\W_{\cL}(Q_0;\F)$.
\item \emph{Periodicity away from a finite set} (\S\ref{sp:sec:automatic}). Cut the positions into blocks of equal length. On an initial segment of each block, the coordinates of the elements of $\W_{\cL}(Q;\F)$ are given by a prefix of one infinite word over a finite alphabet. This word is automatic in two multiplicatively independent bases, hence eventually periodic by Cobham's theorem (Lemma~\ref{lem:periodic_zone}). The points of the circle that lie in no such initial segment, for any block length, form a closed set $E$ that is invariant under multiplication by both bases, hence a finite set of rational points by Furstenberg's theorem (Lemma~\ref{lem:E_rational}). Lemma~\ref{lem:rational_concentration} combines the two by compactness: there is a fixed $P$ such that $x(j+P)=x(j)$ unless $j/Q$ is close to $E$.
\item \emph{Detection} (\S\ref{sp:sec:grid}). 
Let $G$ be a common denominator of the points of $E$ that is divisible by $N\prod_{p\in\calP}p$.
Lemma~\ref{lem:detection} gives, through each position close to a multiple of $Q/G$, a complete row that is either legal or of the form $C^Q_{p,u}$ with $p\in\calP$, in such a way that these rows form a unitriangular matrix. The legal rows remain legal under small shifts, while the rows $C^Q_{p,u}$ with $p\in\calP$ and all their shifts annihilate $x\in\W_{\cL}(Q;\F)$ because $\calP\subseteq\cL$.
Consequently a difference $x(j+h)-x(j)$ with $h$ small that is supported close to the multiples of $Q/G$ must vanish (Corollary~\ref{cor:detection}).
This is the only step using \eqref{e:vanishing_P}; the rows $C^Q_{p,u}$ with $p\in\calP$ handle the denominators for which a legal row is unavailable.
\item \emph{Conclusion} (\S\ref{sp:sec:closure}). After passing to a suitable multiple of the length, the differences $x(j+P)-x(j)$ are supported close to the multiples of $Q/G$ by step (2), so they vanish by step (3). So $x$ has a period shorter than its zero interval (Lemma~\ref{lem:zeros}), and $x=0$.
\end{enumerate}

\subsection{Restrictions and traces}\label{sp:sec:algebra}
Let $Q$ be a positive integer. 
Recall that $C_{q, u}^Q = \sum_{a\bmod{q}} e_{u+a\frac{Q}{q}} = \sum_{d\mid q} R_{d, u}^Q$ for $q\mid Q$ and $u\in\Z$.
By M\"obius inversion, we have $R_{q, u}^Q = \sum_{d\mid q}\mu(q/d) C_{d, u}^Q$, and $\ell_N(d) \leq \ell_N(q)$ for $d\mid q$, so $R_{d, u}^Q$ and $C_{d, u}^Q$ are legal whenever $R_{q, u}^Q$ is.
Hence the legal rows and the legal complete rows span the same subspace of $\F^{\Z/Q\Z}$, for any field $\F$. 

For $m\in \N$ and $0 \leq t < m$, we define the restriction map
\begin{equation}\label{e:eta_def}
\Res_{m, t}^{m Q} : \F^{\Z/(mQ)\Z} \to \F^{\Z/Q\Z}, \qquad (\Res_{m, t}^{m Q} x)(n) = x(mn+t) \quad n\bmod{Q}, \quad x\in \F^{\Z/(mQ)\Z}.
\end{equation}
Like the trace $\Tr_m^{mQ}$ of \eqref{sp:eq:trace}, it goes from length $mQ$ to length $Q$: for both maps, the superscript denotes the input length.

\begin{lemma}\label{lem:maps}
Let $m\in \N$ and $0\leq t < m$.
\begin{enumerate}
\item If $\ell\mid Q$, we have $\Tr_\ell^Q \Res_{m, t}^{m Q} = \Res_{m, t}^{m Q/\ell} \Tr_{\ell}^{mQ}$
and $\Tr_{\ell}^Q \Tr_{m}^{mQ} = \Tr_{m\ell}^{mQ} = \Tr_{m}^{mQ/\ell} \Tr_{\ell}^{mQ}$.
\item Let $\cL$ be a finite set of primes and $Q$ be divisible by $\prod_{p\in \cL} p$.
Then $\Res_{m, t}^{m Q}$ and $\Tr_m^{mQ}$ map $\W_{\cL}(mQ; \F)$ into $\W_{\cL}(Q; \F)$.
\end{enumerate}
\end{lemma}

\begin{proof}
(1) is immediate from the definitions.
By (1), both maps send $\ker\Tr_\ell^{mQ}$ into $\ker\Tr_\ell^{Q}$ for $\ell\in\cL$, so for (2) it suffices to show that they map $\K(mQ;\F)$ into $\K(Q;\F)$.
For $x\in\F^{\Z/(mQ)\Z}$ we have
\[\left<C_{q, u}^Q, \Res_{m, t}^{m Q} x \right> = \left<C_{q, mu+t}^{mQ}, x\right> \qquad\text{and}\qquad \left<C_{q, u}^{Q}, \Tr_m^{mQ} x \right> = \left<C_{mq, u}^{mQ}, x\right>,\]
and if $C_{q, u}^Q$ is legal, then so are $C_{q, mu+t}^{mQ}$ and $C_{mq, u}^{mQ}$: the first because $-\frac{mQ}{\sqrt{\ell_N(q)}}<mu\leq mu+t\leq m(u+1)-1<\frac{mQ}{\sqrt{\ell_N(q)}}-1$, and the second because $\ell_N(mq)$ divides $m^2\ell_N(q)$.
\end{proof} 

\begin{corollary}\label{cor:kernel_growingQ}
Let $\cL$ be a finite set of primes and $Q$ be divisible by $\prod_{p\in \cL}p$. 
If $\W_{\cL}(Q; \F)=0$, then $\W_{\cL}(mQ; \F)=0$ for every positive integer $m$. 
\end{corollary} 
\begin{proof} 
Let $x\in \W_{\cL}(mQ; \F)$. By Lemma \ref{lem:maps}, we have $\Res_{m, t}^{m Q} x\in \W_{\cL}(Q; \F)= 0$ for $0\leq t < m$, that is, $x(mn+t) = 0$ for all $n\bmod{Q}$ and $0\leq t < m$. So $x=0$. 
\end{proof} 

The next lemma allows us to move between $\Q$ and finite fields.

\begin{lemma}\label{sp:lem:field}
Let $\cL$ be a finite set of primes, and let $Q$ be divisible by all elements of $\cL$.
If $\W_{\cL}(Q;\F)=0$ for some field $\F$, then $\W_{\cL}(Q;\Q)=0$.
Conversely, if $\W_{\cL}(Q;\Q)=0$, then $\W_{\cL}(Q;\F)=0$ for every field $\F$ whose characteristic lies outside a finite set of primes depending on $\cL$ and $Q$.
\end{lemma}
\begin{proof}
Let $Z_{\cL}^Q$ be the integer matrix with $Q$ columns whose rows are the legal rows $R_{q,u}^Q$ together with the complete rows $C_{\ell,n}^Q$ for $\ell\in\cL$ and $n\bmod Q/\ell$.
Since $(\Tr_\ell^Qx)(n)=\langle C_{\ell,n}^Q,x\rangle$, the space $\W_{\cL}(Q;\F)$ is the kernel of $Z_{\cL}^Q$ over $\F$, so it is zero if and only if some $Q\times Q$ minor of $Z_{\cL}^Q$ is nonzero in $\F$.
Such a minor is the image in $\F$ of an integer $d$, and it is nonzero in $\F$ if and only if $d\neq 0$ and the characteristic of $\F$ does not divide $d$.
\end{proof}

Finally we record an initial zero interval for $x$.

\begin{lemma}\label{lem:zeros}
Let $x\in\K(Q;\F)$. Then $x(u)=0$ for every integer $u$ with $-\frac{Q}{\sqrt N}<u<\frac{Q}{\sqrt N}-1$. In particular, if $Q/\sqrt N$ is not an integer, then
\begin{equation}\label{sp:eq:zeros}
 x(u)=0\quad\text{for }-\lfloor Q/\sqrt N\rfloor\le u\le \lfloor Q/\sqrt N\rfloor-1.
\end{equation}
\end{lemma}
\begin{proof}
For such $u$ the row $R^Q_{1,u}=e_u$ is legal, because $\ell_N(1)=N$, and $\langle e_u,x\rangle=x(u)$.
\end{proof}

\subsection{Proof of Proposition \ref{sp:prop:initial}}\label{sp:sec:removal}

For a positive integer $Q$ and $m\bmod{Q}$, we define the additive character $\psi_m$ on $\Z/Q\Z$ by
\begin{equation}\label{e:psim_def}
\psi_m(n) = e^{2\pi i nm/Q} \qquad \text{ for } n\bmod{Q}.
\end{equation}

\begin{proof}
We work over $\C$, with $Q=D$. For $p\in \cL$,
\[(\Tr_p^D \psi_m)(n) = \sum_{j\bmod{p}} \psi_m(n+jD/p) = \psi_m(n) \sum_{j\bmod{p}} e^{2\pi i mj/p}
= \begin{cases} p\psi_m(n) & \text{ if } p\mid m, \\ 0 & \text{ otherwise. }\end{cases}\]
So if $x = \sum_{m\bmod{D}} \hat{x}(m) \psi_m$ is the Fourier expansion of $x\in \C^{\Z/D\Z}$, then
\[\Tr_p^D x= p \sum_{m\bmod{D/p}}\hat{x}(mp) \psi_{mp}.\]
Since the additive characters $\psi_{mp}$, $m\bmod{\frac{D}{p}}$, are linearly independent as functions on $\Z/(D/p)\Z$, we have $x\in \ker(\Tr_p^D)$ if and only if $\hat{x}(m)=0$ for every $m$ divisible by $p$.

Take any $x\in \W_{\cL}(D; \C)$.
Then $\hat{x}$ is supported on the residues $m\bmod{D}$ coprime to $D$.
Since $\gcd(D,N)=1$ and $N>1$, the number $D/\sqrt N$ is not an integer, so by Lemma~\ref{lem:zeros} we have $x(j-\lfloor D/\sqrt N\rfloor)=0$ for $0 \leq j \leq 2\lfloor D/\sqrt N\rfloor-1$, that is,
\begin{equation}\label{e:system_initial_vanishing}
\sum_{\substack{m\bmod{D}\\ \gcd(D, m)=1}}  e^{2\pi i mj/D} \cdot \hat{x}(m)\, e^{-2\pi i m\lfloor D/\sqrt N\rfloor/D} = 0 \qquad (0 \leq j \leq 2\lfloor D/\sqrt N\rfloor-1).
\end{equation}
The number of unknowns $\hat{x}(m)$ is $\prod_{p\in\cL}(p-1) = D\rho(\cL)$.
By our assumption, $D\rho(\cL) < 2\frac{D}{\sqrt{N}}$.
Since $\rho(\cL)<\frac{2}{\sqrt N}<1$, the set $\cL$ is not empty, and its elements are odd; so $D\rho(\cL)=\prod_{p\in\cL}(p-1)$ is an even integer.
So we have $D\rho(\cL) \leq 2 \lfloor D/\sqrt N\rfloor$.
The equations \eqref{e:system_initial_vanishing} with $0\leq j < D\rho(\cL)$ form a square system whose matrix $\bigl(e^{2\pi i mj/D}\bigr)_{j, m}$ is a Vandermonde matrix with distinct nodes $e^{2\pi i m/D}$, hence invertible.
So we have $\hat{x}(m)=0$ for all $m$.
Therefore $x=0$, so $\W_{\cL}(D; \C)=0$, and hence $\W_{\cL}(D; \Q)=0$ by Lemma~\ref{sp:lem:field}.
\end{proof}

\subsection{Dimension stabilization}\label{s:finitefield}\label{sp:sec:stationary}

This subsection is the only place where the hypothesis that $\cL\cup\{c_1\}$ and $\cL\cup\{c_2\}$ are vanishing is used, and it is used only through Lemma~\ref{lem:trace_injective}: over a suitable finite field, the traces $\Tr_c^{cQ}$ with $c\in\{c_1,c_2\}$ are injective on $\W_{\cL}$.
It follows that the dimensions of these spaces stabilize, so that at the lengths $Q=Q_0c_1^ac_2^b$ an element of $\W_{\cL}$ is determined by its trace in a fixed finite space $V$ (Lemma~\ref{lem:stabilization}).
Lemma~\ref{lem:Fn} is then a formal consequence of Lemma~\ref{lem:maps}: the restrictions to residue classes are given by fixed linear maps of $V$, and the lemma records how these maps compose.

From now on we continue under the assumptions of Proposition~\ref{sp:prop:two}.
Since $\cL\cup\{c_1\}$ and $\cL\cup\{c_2\}$ are vanishing, by Corollary \ref{cor:kernel_growingQ} there is a positive integer $B$ divisible by $c_1c_2\prod_{p\in \cL}p$ such that
\[\W_{\cL\cup\{c_1\}}(B; \Q) = \W_{\cL\cup\{c_2\}}(B; \Q)=0.\]
We now pass to a finite field so that the stabilized coordinate functionals form a finite alphabet, allowing us to apply Cobham's theorem. By Lemma~\ref{sp:lem:field}, there exists a finite field $\F$ such that
\[\W_{\cL\cup\{c_1\}}(B; \F) = \W_{\cL\cup\{c_2\}}(B; \F)=0.\]
Again by Corollary \ref{cor:kernel_growingQ}, for any $Q$ divisible by $B$, we have
\begin{equation}\label{e:F_fix}
\W_{\cL\cup\{c_1\}}(Q; \F) = \W_{\cL\cup\{c_2\}}(Q; \F)=0.
\end{equation}
We fix $\cL$, $c_1$, $c_2$, $B$ and $\F$, and for $Q$ divisible by every element of $\cL$ we write $\W_Q\coloneq\W_{\cL}(Q;\F)$. 

\begin{lemma}\label{lem:trace_injective} 
For $c\in \{c_1, c_2\}$ and any positive integer $Q$ divisible by $B$, the trace $\Tr_c^{cQ} : \W_{cQ} \to \W_Q$ is injective. 
\end{lemma} 

\begin{proof} 
By Lemma \ref{lem:maps}, we have $\Tr_c^{cQ}\left(\W_{cQ}\right) \subseteq \W_Q$. 
Then 
\[\ker(\Tr_c^{cQ}) \cap \W_{cQ} = \W_{\cL\cup\{c\}}(cQ; \F)=0, \]
by \eqref{e:F_fix}, since $cQ$ is divisible by $B$. 
\end{proof} 

Let $\cI(c_1, c_2) \coloneq \left\{c_1^ac_2^b:\, a, b\in \Z_{\geq 0}\right\}$. 

\begin{lemma}\label{lem:stabilization}
There exist $a_0, b_0\geq 0$ such that, for $Q_0 = Bc_1^{a_0}c_2^{b_0}$, the map
\begin{equation}\label{e:T_a0b0}
\Tr_c^{cQ}: \W_{cQ}\to \W_Q \qquad\text{ for } c\in \{c_1, c_2\} \text{ and all } Q=Q_0m, \, m\in \cI(c_1, c_2)
\end{equation}
is an isomorphism.
\end{lemma}

\begin{proof}
By Lemma \ref{lem:trace_injective}, we have 
\[\dim\W_{Bm_1} \geq \dim\W_{Bm_2} \quad \text{ whenever } m_1\mid m_2, \quad m_1, m_2\in \cI(c_1, c_2).\]
Choose $a_0,b_0\geq0$ for which $\dim\W_{Bc_1^{a_0}c_2^{b_0}}$ is minimal, and let $Q_0 = Bc_1^{a_0}c_2^{b_0}$.
Then $\W_{Q_0m}$ has the same dimension for every $m\in \cI(c_1, c_2)$.
So for $Q=Q_0m$ and $c\in\{c_1,c_2\}$, the map $\Tr_c^{cQ}$ of \eqref{e:T_a0b0} is an injective linear map between spaces of the same finite dimension, hence an isomorphism.
\end{proof}

We fix $Q_0$ as in Lemma~\ref{lem:stabilization} and let $V\coloneq \W_{Q_0}$.

For $m\in \cI(c_1, c_2)$ and $Q=Q_0m$, the map
\[\Tr_m^Q: \W_Q\to V = \W_{Q_0}\]
is a composition of isomorphisms \eqref{e:T_a0b0}, by Lemma \ref{lem:maps}(1).
Let
\[I_Q : V\to \W_Q\]
be its inverse.
For $c\in \{c_1, c_2\}$, again by Lemma \ref{lem:maps}(1), we have
\begin{equation}\label{e:TmcQ_split}
\Tr_{cm}^{cQ} = \Tr_m^Q \Tr_c^{cQ}.
\end{equation}
Then $\Tr_{m}^Q (\Tr_c^{cQ} I_{cQ}) = \Tr_{cm}^{cQ} I_{cQ}$ is the identity map, so we have
\begin{equation}\label{e:IQ_split}
\Tr_c^{cQ} I_{cQ} = I_Q.
\end{equation}
Iterating \eqref{e:IQ_split} gives $\Tr_n^{nQ}I_{nQ}=I_Q$ for every $n\in\cI(c_1,c_2)$.

For $n\in \cI(c_1, c_2)$ and $0\leq i < n$, define
\[F_Q^{(n)}(i) = \Tr_{m}^Q \circ \Res_{n, i}^{n Q} \circ I_{nQ}\in \End_{\F}(V).\]
By Lemma \ref{lem:maps}, $\Res_{n, i}^{n Q} \W_{nQ} \subseteq \W_Q$, so $F_Q^{(n)}(i)$ is well defined on $V$.
The finite ring $\End_{\F}(V)$ acts on the right on the finite set $V^*=\Hom_{\F}(V,\F)$ by $\alpha A = \alpha\circ A$.

We define the coordinate word of length $Q$ with letters in $V^*$ by
\[w_Q(j)\colon V\to\F,\qquad \xi\mapsto (I_Q\xi)(j) \qquad (0\leq j < Q). \]
Thus $w_Q(j)$ is the $j$-th coordinate of an element of $\W_Q$ as a linear form in its trace; it is zero if and only if every such element vanishes at $j$.

\begin{lemma}\label{lem:Fn}
For $n\in \cI(c_1, c_2)$ and $0\leq i < n$, the endomorphism $F_Q^{(n)}(i)$ does not depend on $Q$; we write $F^{(n)}(i) = F_Q^{(n)}(i)$. Moreover, the following hold.
\begin{enumerate} 
\item For $n_1, n_2 \in \cI(c_1, c_2)$, $0\leq i_1 < n_1$ and $0\leq i_2 < n_2$, 
\begin{equation}\label{e:Fn_add}
F^{(n_1n_2)} (n_2i_1+i_2) = F^{(n_1)}(i_1)F^{(n_2)}(i_2). 
\end{equation} 
\item Recall that $Q=Q_0m$ for $m\in \cI(c_1, c_2)$. 
For $n\in \cI(c_1, c_2)$, $0\leq j < Q$ and $0\leq i < n$, 
\begin{equation}\label{e:wnQ_Fn}
w_{nQ}(nj+i) = w_Q(j) F^{(n)}(i). 
\end{equation} 
\item For any $x\in \W_Q$, 
\[x(j) = w_Q(j) (\Tr_m^{Q} x).\]
\end{enumerate} 
\end{lemma} 

\begin{proof} 
We first show that $F_Q^{(n)}(i)$ does not depend on $Q$.
Because $Q = mQ_0$ and $m\in \cI(c_1, c_2)$, it is enough to show $F_{cQ}^{(n)}(i) = F_{Q}^{(n)}(i)$ for $c\in \{c_1, c_2\}$. 
By \eqref{e:TmcQ_split}, 
\[F_{cQ}^{(n)}(i) = \Tr_{cm}^{cQ} \Res_{n, i}^{n cQ} I_{ncQ} = \Tr_{m}^{Q} (\Tr_c^{cQ} \Res_{n, i}^{n cQ}) I_{ncQ}.\]
By Lemma \ref{lem:maps}, 
\[\Tr_c^{cQ} \Res_{n, i}^{n cQ} = \Res_{n, i}^{n Q} \Tr_c^{ncQ}.\]
So by \eqref{e:IQ_split}, 
\[F_{cQ}^{(n)}(i) = \Tr_{m}^{Q} \Res_{n, i}^{n Q} (\Tr_c^{ncQ} I_{ncQ}) = \Tr_{m}^{Q} \Res_{n, i}^{n Q} I_{nQ} = F_Q^{(n)}(i), \]
i.e., $F^{(n)}_Q$ does not depend on $Q$. 

(1) By \eqref{e:IQ_split}, we have $I_{n_2Q} = \Tr_{n_1}^{n_1n_2Q} I_{n_1n_2Q}$. 
So we have 
\[F^{(n_2)}(i_2) = \Tr_{m}^Q \circ \Res_{n_2, i_2}^{n_2 Q} \circ \Tr_{n_1}^{n_1n_2Q} \circ I_{n_1n_2Q}.\]
By Lemma \ref{lem:maps}, $\Res_{n_2, i_2}^{n_2 Q} \Tr_{n_1}^{n_1n_2 Q} = \Tr_{n_1}^{n_1Q} \Res_{n_2, i_2}^{n_2 n_1Q}$, so we have 
\[F^{(n_2)}(i_2)= \Tr_m^Q \circ \Tr_{n_1}^{n_1Q} \circ \Res_{n_2, i_2}^{n_2 n_1Q} \circ I_{n_1n_2Q} = \Tr_{mn_1}^{n_1Q} \circ \Res_{n_2, i_2}^{n_2 n_1Q} \circ I_{n_1n_2Q}.\]
Note that $I_{n_1Q} \Tr_{mn_1}^{n_1Q}$ is the identity map on $\W_{n_1Q}$. 
So 
\[F^{(n_1)}(i_1) F^{(n_2)}(i_2) = \Tr_{m}^Q \circ \Res_{n_1, i_1}^{n_1 Q} \circ \Res_{n_2, i_2}^{n_2 n_1Q} \circ I_{n_1n_2Q}.\]
Since $\Res_{n_1, i_1}^{n_1 Q}\Res_{n_2, i_2}^{n_2 n_1Q} = \Res_{n_1n_2, n_2i_1+i_2}^{n_1n_2 Q}$ (both sides send $x$ to $j\mapsto x(n_1n_2 j + n_2 i_1+ i_2)$), we have 
\[F^{(n_1)}(i_1) F^{(n_2)}(i_2) = \Tr_{m}^Q \circ \Res_{n_1n_2, n_2i_1+i_2}^{n_1n_2 Q} \circ I_{n_1n_2Q} = F^{(n_1n_2)}(n_2i_1+i_2).\]

(2) Let $\xi\in V$. Since $x(nj+i) = (\Res_{n, i}^{n Q} x)(j)$ for $x\in\F^{\Z/(nQ)\Z}$, we have 
\[w_{nQ}(nj+i)(\xi) = \bigl(\Res_{n, i}^{n Q} I_{nQ}\xi\bigr)(j).\]
Since $\Res_{n, i}^{n Q} I_{nQ}\xi\in\W_Q$ and $I_Q \Tr_m^Q$ is the identity on $\W_Q$, we have 
\[w_{nQ}(nj+i)(\xi) = \bigl(I_Q\Tr_m^Q \Res_{n, i}^{n Q} I_{nQ}\xi\bigr)(j) = w_Q(j)\bigl(F^{(n)}(i)\xi\bigr),\]
that is, $w_{nQ}(nj+i) = w_Q(j) F^{(n)}(i)$.

(3) Since $\Tr_m^Q$ is an isomorphism and $I_Q$ is its inverse, we have 
\[w_Q(j)(\Tr_m^Q x) = (I_Q\Tr_m^Q x)(j) = x(j).\]
\end{proof} 

\subsection{Periodicity away from a finite set}\label{s:automatic}\label{sp:sec:automatic}

The goal of this subsection is to prove Lemma~\ref{lem:rational_concentration}. This lemma records that $x(j + P) = x(j)$ for a fixed $P$ and $x \in \W_Q$ unless $j/Q$ lies near a finite subset $E$ of the circle.
Its proof has two independent parts.
Lemma~\ref{lem:periodic_zone} uses Cobham's theorem to show that $x(j+P)=x(j)$ when $j$ lies in an initial segment of a block of consecutive positions, for a whole family of block lengths.
Lemma~\ref{lem:E_rational} uses Furstenberg's theorem to show that the points of the circle covered by no such segment form a finite set $E$ of rational points.
The lemma then follows from Lemma~\ref{lem:periodic_zone} by compactness; its proof does not use that $E$ is finite, but without Lemma~\ref{lem:E_rational} the statement could be empty.

Let $q\geq 2$. For an integer $n\geq 0$ let $[n]_q = (r_1, \ldots, r_T)$ for $n=\sum_{t=1}^T r_t q^{T-t}$ denote the base $q$ expansion of $n$.

Let $A$ be a finite set. 
A sequence $u: \Z_{\geq 0} \to A$ is $q$-automatic if there is a finite set $S$ of states, an initial state $s_0\in S$, a transition map $S\times \{0, \ldots, q-1\}\to S$, $(s, r)\mapsto s\cdot r$, 
and an output map $\Omega: S\to A$ satisfying the following:  whenever $[n]_q=(r_1, \ldots, r_T)$, 
\[u(n) = \Omega\left((\cdots ((s_0\cdot r_1)\cdot r_2)\cdots )\cdot r_T\right).\]
A sequence $u(0), u(1), u(2), \ldots$ is eventually periodic if there are integers $i_0\geq 0$ and $P\geq 1$ with 
\[u(i+P) = u(i)\qquad\text{ for all } i \geq i_0.\]
Two integers $q_1, q_2\geq 2$ are multiplicatively independent if $q_1^{a_1}\neq q_2^{a_2}$ for all positive integers $a_1, a_2$. 

The following is Cobham's theorem \cite{Cobham}, proved there for sets of integers. 
The version for sequences below follows by applying it to each set $\{n: \, u(n)=a\}$ (see also \cite[Chapter 11]{AS03}). 
\begin{theorem}\label{thm:Cobham}
Let $q_1, q_2\geq 2$ be multiplicatively independent. 
A sequence that is both $q_1$-automatic and $q_2$-automatic is eventually periodic. 
\end{theorem}
The following is a consequence of \cite[Lemma IV.2 and Theorem IV.1]{Furstenberg}. 
\begin{theorem}\label{thm:Furstenberg}
Let $q_1, q_2\geq 2$ be multiplicatively independent. 
A closed subset $A$ of $\R/\Z$ with $q_1 A\subseteq A$ and $q_2 A\subseteq A$ is either a finite subset of $\Q/\Z$ or equal to $\R/\Z$. 
\end{theorem}

\begin{proof}
The semigroup $\Sigma=\{q_1^aq_2^b:\, a, b\geq 0\}$ is non-lacunary in the sense of \cite[Definition IV.1]{Furstenberg}: if all its members were powers of a single integer $g$, say $q_1=g^i$ and $q_2=g^j$, then $q_1^j=q_2^i$.
Suppose that $A\neq \R/\Z$.
If $A$ contained an irrational point $y$, then $\Sigma y\subseteq A$ would be dense in $\R/\Z$ by \cite[Theorem IV.1]{Furstenberg}, and the closed set $A$ would be $\R/\Z$.
Hence $A\subseteq \Q/\Z$.
If $A$ were infinite, it would have a non-isolated point $y$, and $hy=0$ for some integer $h\geq 1$.
The set $hA$ is closed and $\Sigma$-invariant, and $0$ is a non-isolated point of it.
So $hA=\R/\Z$ by \cite[Lemma IV.2]{Furstenberg}, which is impossible because $hA\subseteq \Q/\Z$ is countable.
Hence $A$ is a finite subset of $\Q/\Z$.
\end{proof}

When the length is multiplied by $c\in\{c_1,c_2\}$, each letter $w_Q(j)$ is replaced by a block of $c$ letters whose first letter is $w_{cQ}(cj)=w_Q(j)F^{(c)}(0)$, in general different from $w_Q(j)$.
Passing to powers $q_1$, $q_2$ of $c_1$, $c_2$ and to letters in a set $V^*\Pi$ removes this, so that the blocks begin with longer and longer prefixes of one infinite word.
Since $\End_{\F}(V)$ is finite, for every element $e\in \End_{\F}(V)$, there exists $s\in \N$ such that $e^s$ is idempotent, i.e., $e^s e^s = e^s$. 
Choose $s_1, s_2\in \N$ such that $F^{(c_1)}(0)^{s_1}$ and $F^{(c_2)}(0)^{s_2}$ are idempotents. 
Let 
\[q_1= c_1^{s_1}, \qquad q_2 = c_2^{s_2}, \qquad \cI(q_1, q_2) = \{q_1^{a}q_2^b:\, a, b\in \Z_{\geq 0}\} \subseteq \cI(c_1, c_2).\]
The integers $q_1$ and $q_2$ are multiplicatively independent since $c_1\neq c_2$ are prime numbers. 
By \eqref{e:Fn_add}, taking $i_1=i_2=0$, we have 
\[F^{(n_1n_2)}(0) = F^{(n_1)}(0) F^{(n_2)}(0)\qquad n_1, n_2\in \cI(c_1, c_2).\]
So, for $c\in \{c_1, c_2\}$ and $s\in \Z_{\geq 0}$, we have 
\[F^{(c^s)}(0) = (F^{(c)}(0))^s.\]
Let 
\[\Pi_1 = F^{(q_1)}(0) = (F^{(c_1)}(0))^{s_1}, \qquad \Pi_2 = F^{(q_2)}(0) = (F^{(c_2)}(0))^{s_2}.\]
By the choice of $s_1$ and $s_2$, these are idempotents, and they commute:
\[\Pi_1 \Pi_2 = F^{(q_1)}(0)F^{(q_2)}(0) = F^{(q_1q_2)}(0) = \Pi_2\Pi_1.\]
Let $\Pi = \Pi_1\Pi_2$. Then $\Pi^2=\Pi$ and $\Pi \Pi_1 = \Pi\Pi_2 = \Pi$. 
So 
\begin{equation}\label{e:Pi}
\Pi \text{ is idempotent} \qquad\text{ and } \qquad \Pi F^{(n)}(0) = \Pi \qquad\text{ for any } n\in \cI(q_1, q_2).
\end{equation}

Let 
\[V^* \Pi = \left\{\beta \Pi:\, \beta\in V^*\right\}.\]
Then $\alpha\in V^*\Pi$ if and only if $\alpha\Pi = \alpha$. 

\begin{lemma}\label{lem:V*Pi}
Let $\alpha\in V^*\Pi$. 
\begin{enumerate} 
\item For $i\geq 0$, the letter $u_\alpha(i) = \alpha F^{(n)}(i)$  does not depend on the choice of $n\in \cI(q_1, q_2)$ with $n>i$, i.e., $u_\alpha(i)$ is well-defined. 

\item $u_\alpha(0)=\alpha$ and 
\begin{equation}\label{e:ualpha_F}
u_{\alpha}(qi +r) = u_\alpha (i) F^{(q)}(r)\qquad\text{ for }q\in \{q_1, q_2\}, \, i\geq 0\text{ and } 0 \leq r < q.
\end{equation} 
\end{enumerate}
\end{lemma} 

\begin{proof} 
(1) Let $n, n'\in \cI(q_1, q_2)$ with $n, n'>i$.
By \eqref{e:Fn_add},
\[F^{(nn')}(i) = F^{(n)}(0)F^{(n')}(i) = F^{(n')}(0) F^{(n)}(i).\]
Since $\alpha = \alpha\Pi$ and $\Pi F^{(n)}(0) = \Pi = \Pi F^{(n')}(0)$ by \eqref{e:Pi}, applying $\alpha$ to both sides gives $\alpha F^{(n')}(i) = \alpha F^{(n)}(i)$.

(2) For $n\in \cI(q_1, q_2)$, 
\[u_\alpha(0) = \alpha F^{(n)}(0) = \alpha \Pi F^{(n)}(0) = \alpha \Pi = \alpha.\]
By \eqref{e:Fn_add}, for $n\in \cI(q_1, q_2)$ with $n>i$, 
\[u_\alpha(qi+r) = \alpha F^{(nq)}(qi+r) = \alpha F^{(n)}(i) F^{(q)}(r) = u_\alpha(i) F^{(q)}(r).\]
\end{proof} 

\begin{lemma}\label{lem:ualpha_even_periodic}
The sequence $u_\alpha$ is $q_1$-automatic and $q_2$-automatic, hence eventually periodic. 
There exist $i_0\in \Z_{\geq 0}$ and $P\in \N$ such that 
\begin{equation}\label{e:ualpha_even_periodic}
u_\alpha(i+P) = u_\alpha(i) \qquad\text{ for all } \alpha\in V^*\Pi \text{ and for all } i \geq i_0.
\end{equation} 
\end{lemma} 

\begin{proof} 
Let $q\in \{q_1, q_2\}$. 
Take $V^*$ as the set of states, $\alpha$ as the initial state, the transition map
\[\beta\cdot r = \beta F^{(q)}(r),  \qquad \beta\in V^*, \quad r\in \{0, 1, \ldots, q-1\},\]
and the identity map of $V^*$ as the output map.
We claim that the state reached after reading the digits $r_1, \ldots, r_T$ is $u_\alpha(n)$, where $n=\sum_{t=1}^T r_t q^{T-t}$.
For $T=0$ this is $u_\alpha(0) = \alpha$.
If it holds for $r_1, \ldots, r_T$, then after reading one more digit $r_{T+1}$ the state is $u_\alpha(n) F^{(q)}(r_{T+1}) = u_\alpha(qn+r_{T+1})$ by \eqref{e:ualpha_F}, which proves the claim by induction on $T$.
Thus $u_\alpha$ is $q$-automatic for $q\in \{q_1, q_2\}$, and by Theorem \ref{thm:Cobham} it is eventually periodic.
Since $V^*\Pi$ is finite, we can choose $i_0$ and $P$ that work for all $\alpha\in V^*\Pi$ simultaneously, which is \eqref{e:ualpha_even_periodic}.
\end{proof} 

Fix $P$ and $i_0$ as in Lemma \ref{lem:ualpha_even_periodic}. 
For $x\in \F^{\Z/Q\Z}$, let 
\begin{equation}\label{e:DeltaP_def}
(\Delta_P x )(j) = x(j+P)-x(j), \qquad 0 \leq j < Q.
\end{equation}
Differences are natural here because a shift of a legal row is usually legal: if $R^Q_{q,u}$ and $R^Q_{q,u+P}$ are both legal, then $\langle R^Q_{q,u},\Delta_Px\rangle=\langle R^Q_{q,u+P},x\rangle-\langle R^Q_{q,u},x\rangle=0$ for every $x\in\K(Q;\F)$. 

\begin{lemma}\label{lem:periodic_zone}
Let $m, n\in \cI(q_1, q_2)$ and $Q=Q_0 mq_1q_2 n$. 
Let $0\leq j < Q$ and let 
\[\theta \coloneq \frac{Q_0 mj}{Q} - \left\lfloor\frac{Q_0mj}{Q}\right\rfloor \in [0, 1), \]
be the fractional part of $\frac{Q_0 mj}{Q}$. 
If 
\begin{equation}\label{e:periodic_zone}
\theta\in \left[\frac{i_0}{q_1q_2 n} , \frac{1}{q_1q_2} - \frac{P}{q_1q_2n} \right), 
\end{equation} 
then $w_Q(j+P)= w_Q(j)$, and therefore $(\Delta_Px)(j) =0$ for every $x\in \W_Q$. 
\end{lemma} 

In other words, cut the positions $0,\ldots,Q-1$ into $Q_0m$ blocks of $q_1q_2n$ consecutive positions; condition \eqref{e:periodic_zone} says that the offset $i$ of $j$ in its block satisfies $i_0\leq i<n-P$, so that $j$ and $j+P$ are both among the first $n$ positions of the block.

\begin{proof} 
Let $Q_1=Q_0m$ and $M=q_1q_2 n$, so that $Q = Q_1 M$.
Write $j= j_1 M + i$ with $0 \leq j_1 < Q_1$ and $0\leq i < M$. 
Then 
\[\frac{Q_0 mj}{Q} = \frac{Q_1(j_1M+i)}{Q_1M} = j_1 + \frac{i}{M}\qquad \text{ so } \theta=\frac{i}{M}. \]
We assume \eqref{e:periodic_zone}, so $i \in [i_0, n-P)$. 

Let $\alpha = w_{Q_1}(j_1) \Pi\in V^*\Pi$. 
For every $0\leq i' < n$, by \eqref{e:wnQ_Fn}, \eqref{e:Fn_add} and Lemma \ref{lem:V*Pi}, 
\begin{align*}
w_Q(j_1 M + i') & = w_{Q_1}(j_1) F^{(M)}(i') = w_{Q_1}(j_1) F^{(q_1q_2)}(0)F^{(n)}(i') = w_{Q_1}(j_1)\Pi F^{(n)}(i') 
\\ & = \alpha F^{(n)}(i') = u_{\alpha}(i').
\end{align*}
By Lemma \ref{lem:ualpha_even_periodic}, we have 
\[w_Q(j_1 M + i' + P) = u_{\alpha}(i'+P) = u_{\alpha}(i') = w_{Q}(j_1M + i')\qquad \text{ for all } i_0 \leq i' < n-P. \]
Taking $i'=i$, we have $i'\geq i_0$ and $i'+P< n$, so 
\[w_Q(j) = w_Q(j_1M+i) = w_Q(j_1M+i+P) = w_Q(j+P).\]
By Lemma \ref{lem:Fn} (3), for any $x\in \W_Q$, we have 
\[x(j) = w_Q(j) (\Tr_{Q/Q_0}^Qx) = w_Q(j+P)(\Tr_{Q/Q_0}^Qx) = x(j+P), \]
so $(\Delta_Px)(j)=0$. 
\end{proof} 

For $t_1, t_2\in \R/\Z$, let $\dist(t_1, t_2) = \min_{j\in \Z}  |t_1-t_2-j|$ be the distance on the circle, and for a finite non-empty set $E\subset \R/\Z$, let $\dist(t, E) = \min_{e\in E} \dist(t, e)$. 

Let $\mfm = (0, 1/(q_1q_2)) \subset \R/\Z$ be an open interval and  
\begin{equation}\label{e:E_def}
E = \left\{t\in \R/\Z:\, Q_0m t\not\in \mfm \text{ for all } m\in \cI(q_1, q_2)\right\}. 
\end{equation}
For a position $j$ at length $Q$, the number $\theta$ of Lemma~\ref{lem:periodic_zone} is the fractional part of $Q_0mt$ with $t=j/Q$. So, up to the two ends of the interval \eqref{e:periodic_zone}, that lemma applies to $j$ for some $m$ unless the point $j/Q$ of $\R/\Z$ lies in $E$; see Figure~\ref{fig:zones}. 
\begin{figure}[ht]
\centering
\begin{tikzpicture}[x=12.4cm,y=1cm]
\foreach \m/\y in {1/0, 2/-0.5, 3/-1.0, 4/-1.5, 6/-2.0, 8/-2.5, 9/-3.0, 12/-3.5} {
  \node[left] at (0,\y+0.12) {\footnotesize $m=\m$};
  \pgfmathtruncatemacro{\last}{\m-1}
  \foreach \k in {0,...,\last} {
    \fill[black!35] (\k/\m,\y) rectangle (\k/\m+1/6/\m,\y+0.24);
    \draw[line width=0.3pt] (\k/\m,\y) -- (\k/\m,\y+0.24);
  }
  \draw[line width=0.3pt] (0,\y) -- (1,\y);
  \draw[line width=0.3pt] (1,\y) -- (1,\y+0.24);
}
\node[left] at (0,-4.2+0.12) {\footnotesize $m\leq 12$};
\foreach \m in {1, 2, 3, 4, 6, 8, 9, 12} {
  \pgfmathtruncatemacro{\last}{\m-1}
  \foreach \k in {0,...,\last} {
    \fill[black!35] (\k/\m,-4.2) rectangle (\k/\m+1/6/\m,-4.2+0.24);
  }
}
\draw[line width=0.3pt] (0,-4.2) -- (1,-4.2);
\node[left] at (0,-4.9+0.12) {\footnotesize all $m$};
\fill[black!35] (0,-4.9) rectangle (1,-4.9+0.24);
\draw[line width=0.3pt] (0,-4.9) -- (1,-4.9);
\foreach \t in {0, 1/6, 1/5, 5/24, 5/23, 1/4, 7/24, 7/23, 5/16, 1/3, 2/5, 5/12, 10/23, 7/16, 7/15, 11/23, 23/48, 1/2, 7/12, 3/5, 29/48, 14/23, 5/8, 31/48, 15/23, 2/3, 11/15, 17/23, 3/4, 4/5, 13/16, 19/23, 5/6, 13/15, 20/23, 7/8, 21/23, 11/12, 14/15, 15/16, 22/23, 23/24, 47/48, 1} {
  \draw[line width=0.5pt] (\t,-4.9-0.08) -- (\t,-4.9+0.32);
}
\foreach \t/\l in {0/0, 0.25/\frac14, 0.5/\frac12, 0.75/\frac34, 1/1} {
  \node[below] at (\t,-5.0) {\footnotesize $\l$};
}
\end{tikzpicture}
\caption{The set $E$ of \eqref{e:E_def} for $Q_0=1$, $q_1=2$, $q_2=3$.
Each of the first eight rows shows the circle $\R/\Z$ cut into $m$ arcs, for one $m\in\cI(2,3)$, with the first sixth of each arc shaded: these are the points $t$ with $mt\in\mfm$, to which Lemma~\ref{lem:periodic_zone} applies.
The row ``$m\leq12$'' is the union of these eight rows.
In the last row the union is taken over all $m\in\cI(2,3)$: every point is shaded except the $43$ rational points marked, which form $E$.}
\label{fig:zones}
\end{figure}

\begin{lemma}\label{lem:E_rational}
The set $E$ is a finite subset of $\Q/\Z$ that contains $0$. 
\end{lemma} 

\begin{proof}
Clearly $0\in E$, and $E\neq \R/\Z$ because $Q_0t\in \mfm$ when $0<t<\frac{1}{Q_0q_1q_2}$.
For each $m\in \cI(q_1, q_2)$, the set of $t\in \R/\Z$ with $Q_0mt\in \mfm$ is open, so $E$ is closed.
If $t\in E$ and $q\in \{q_1, q_2\}$, then $Q_0 m (qt) = Q_0(mq)t \not\in \mfm$ for every $m\in \cI(q_1, q_2)$, since $mq\in \cI(q_1, q_2)$.
So $q_1E\subseteq E$ and $q_2 E\subseteq E$.
Since $q_1, q_2$ are multiplicatively independent, $E$ is a finite subset of $\Q/\Z$ by Theorem \ref{thm:Furstenberg}. 

\end{proof} 

\begin{lemma}\label{lem:rational_concentration}
For any $\delta>0$, there is $a_0 \in \Z_{\geq 0}$ such that for any $a, b\geq a_0$ and $Q=Q_0 q_1^{a}q_2^b$, every $x\in \W_Q$ satisfies 
\begin{equation}\label{e:rational_concentration}
x(j+P)=x(j) \qquad\text{ for all } 0\leq j<Q \text{ with } \dist(j/Q, E) \geq \delta. 
\end{equation} 
\end{lemma} 

\begin{proof} 
Let $\delta>0$.
The set $K=\{z\in \R/\Z:\, \dist(z, E) \geq \delta\}$ is closed, hence compact, and it is disjoint from $E$.
For $m\in \cI(q_1, q_2)$ and $\eta\in (0, \frac{1}{2q_1q_2})$, let
\[O_{m, \eta} = \left\{y\in \R/\Z:\, Q_0my +\ell \in (\eta, \frac{1}{q_1q_2}-\eta)\quad\text{ for some } \ell\in \Z\right\}.\]
These sets are open, $O_{m, \eta}$ grows as $\eta$ decreases, and by the definition \eqref{e:E_def} of $E$,
\[\bigcup_{m, \eta} O_{m, \eta} = (\R/\Z)\setminus E \supseteq K.\]
Since $K$ is compact, there exist $m_1, \ldots, m_r\in \cI(q_1, q_2)$ and $\eta\in (0, \frac{1}{2q_1q_2})$ such that
\begin{equation}\label{e:complement_meta}
\bigcup_{i=1}^r O_{m_i, \eta} \supseteq K.
\end{equation} 
Write $m_i = q_1^{a_i}q_2^{b_i}$ for some $a_i, b_i\geq 0$ and $\lambda=\frac{\max\{i_0, P\}}{q_1q_2\eta}$.  
We take $a_0> \max\{a_i, b_i:\, 1\leq i \leq r\}$ such that $2^{a_0-1-a_i} \geq\lambda$ for all $1\leq i \leq r$. 
Take $a, b\geq a_0$ and $Q=Q_0q_1^{a} q_2^{b}$. 
For each $1\leq i \leq r$, we have 
\[Q = Q_0 m_i q_1q_2 n_i \qquad \text{ for } n_i = q_1^{a-a_i-1} q_2^{b-b_i -1} \in \cI(q_1, q_2).\]
Since $q_1\geq 2$, and by assumption, we have $n_i \geq \lambda$. 
Then we have 
\[\left(\eta, \frac{1}{q_1q_2}-\eta\right) \subseteq \left[\frac{i_0}{q_1q_2 n_i}, \frac{1}{q_1q_2}-\frac{P}{q_1q_2n_i}\right).\]
If $0\leq j < Q$, and $\dist(j/Q, E) \geq \delta$, then by \eqref{e:complement_meta}, the fractional part of $\frac{Q_0m_i j}{Q}$ lies in $\left(\eta, \frac{1}{q_1q_2}-\eta\right)$ for some $i$. 
By Lemma \ref{lem:periodic_zone}, we have $(\Delta_Px)(j)=0$ for any $x\in \W_Q$. 
So we have \eqref{e:rational_concentration}.
\end{proof} 

\subsection{Detecting rows near the finite set \texorpdfstring{$E$}{E}}\label{s:detection}\label{sp:sec:grid}

Assume \eqref{e:vanishing_P}, and put
$N_s=\prod_{p\mid s}p^{\ord_p(N)}$.
In Lemma~\ref{lem:detection} below, a legal complete row works at a denominator $s>2$ exactly when $s^2>N_s$. 
At every other denominator $s>2$, \eqref{e:vanishing_P} supplies a prime $p\in\calP$ not dividing $s$. 
The row $C^Q_{p,u}$ through that position has one entry at denominator $s$ and all its other entries at larger denominators, so the restriction matrix in Lemma~\ref{lem:detection} is still unitriangular.

\begin{lemma}\label{lem:detection}
Let $G$ be divisible by $N\prod_{p\in\calP}p$. There exist $H_0\ge1$ and
$0<\varepsilon<1/4$ with the following property. If $H_0\mid v$,
$v\ge1/\varepsilon$, and $Q=Gv$, put
\[\scrE_v=\{av+t\bmod Q:a\bmod G,\ t\in\Z, |t|\le\varepsilon v\}.
\]
For each label $(a,t)$, that is, each pair with $a\bmod G$, $t\in\Z$ and $|t|\le\varepsilon v$, there is a complete row $C^Q_{q,u}$
whose coefficient at $av+t$ is one, such that either
\begin{enumerate}
\item $q\mid Q$ and
\begin{equation}\label{e:margin}
 -Q/\sqrt{\ell_N(q)}<u,\qquad
 u+\varepsilon v<Q/\sqrt{\ell_N(q)}-1;
\end{equation}
\item $q\in\calP$, with no condition on $u$.
\end{enumerate}
The matrix restricting these rows to $\scrE_v$ is
unitriangular after a simultaneous ordering of rows and columns.
In particular its determinant is one over every field.
\end{lemma}

\begin{proof}
For $a\bmod G$ let
\[
 s=\frac G{\gcd(a,G)},\qquad b=\frac a{\gcd(a,G)}.
\]
For $s>1$ this expresses $a/G=b/s$ in lowest terms. We treat
the exceptional denominators $s=1,2$ separately. For $s>2$, call a label good if $s^2>N_s$
and bad otherwise.

For every good label define
\[
 m_s=\prod_{p\mid s}p^{\ord_p(G)},\quad
 d_s=G/m_s,\quad c_s=m_s/s.
\]
Then $G=m_sd_s$, $(m_s,d_s)=1$, and $a=d_sbc_s$.
In particular $(d_sb,s)=1$. Choose a positive integer $r_a$ with
\begin{equation}\label{eq:ra}
 r_a d_s b\equiv1\pmod s,\qquad (r_a,G)=1.
\end{equation}
To do so, prescribe the inverse of $d_sb$ modulo $s$ and prescribe
$r_a\equiv1\pmod p$ for each prime $p\mid G$ not dividing $s$.
The Chinese remainder theorem supplies a positive solution, and
the inverse congruence ensures that it is also a unit at every
prime dividing $s$. Multiplying the first congruence by $c_s$
gives $r_a a\equiv c_s\pmod{m_s}$.

Let $H_0$ be the least common multiple of the finitely many $r_a$,
with value one if there are none, and put
$R=\max(\{1\}\cup\{r_a:a\text{ good}\})$.
Choose $\varepsilon$ so that
\begin{equation}\label{eq:epsilon}
 0<\varepsilon<\frac1{4R},\qquad
 4\varepsilon\le
 \left(\frac s{\sqrt{N_s}}-1\right)\frac{c_s}{r_a}
 \quad\text{for every good label}.
\end{equation}
This is possible because all the finitely many right sides are positive.

For a good label $(a,t)$, choose
\begin{equation}\label{eq:goodrow}
 q=d_s r_a,\qquad u=c_s v/r_a+t.
\end{equation}
Since $r_a\mid v$ and $d_s\mid G$, both the integrality of $u$
and the divisibility $q\mid Q$ hold. The conditions $N\mid G$
and $(r_a,G)=1$ give, prime by prime,
\[
 \ell_N(d_s r_a)=r_a^2d_s^2N_s.
\]
Writing $\alpha=Q/(r_ad_s\sqrt{N_s})$, we have
\[
 \alpha-\frac{c_sv}{r_a}
 =\left(\frac s{\sqrt{N_s}}-1\right)\frac{c_sv}{r_a}
 \ge4\varepsilon v.
\]
Also $u\ge v/R-\varepsilon v>0>-\alpha$, and
\[
 u+\varepsilon v\le c_sv/r_a+2\varepsilon v
 \le\alpha-2\varepsilon v<\alpha-1,
\]
because $\varepsilon v\ge1$. Thus \eqref{e:margin} holds.

For a bad label, \eqref{e:vanishing_P} implies that some
$p=p(a)\in\calP$ does not divide $s$: otherwise every prime of $\calP$ would divide $s$, contradicting that condition. Choose
\begin{equation}\label{eq:badrow}
 q=p(a),\qquad u=av+t.
\end{equation}
This row is of type~(2) in the lemma, and it is defined because $p(a)\mid G\mid Q$.

For $s=1$, so $a=0$, use $C^Q_{1,t}$. For $s=2$, so $a=G/2$,
use $C^Q_{2,t}$. In either case $s\mid G$, and
$\ell_N(s)\mid G^2$, whence $Q/\sqrt{\ell_N(s)}\ge v$.
Consequently
\[
 t\ge-\varepsilon v>-v\ge-Q/\sqrt{\ell_N(s)},
 \qquad t+\varepsilon v\le2\varepsilon v<v-1
 \le Q/\sqrt{\ell_N(s)}-1.
\]
Here $(1-2\varepsilon)v>2$ follows from
$\varepsilon<1/4$ and $v\ge1/\varepsilon$. This proves the
strict margins for the exceptional denominators $s=1,2$.

It remains to prove the matrix assertion. The labels $(a,t)$ are
distinct modulo $Q$, since a coincidence forces $t-t'$ to be a
multiple of $v$, while $|t-t'|\le2\varepsilon v<v$.

For a good row \eqref{eq:goodrow}, its coefficient at $a'v+t'$
is the number of $j\bmod d_s r_a$ satisfying
\[
 \frac v{r_a}(c_s+jm_s)+t\equiv a'v+t'\pmod{Gv}.
\]
Reduction modulo $v/r_a$ forces $t=t'$, because
$|t-t'|\le2\varepsilon v<v/R\le v/r_a$.
The remaining congruence is
\[
 c_s+jm_s\equiv a'r_a\pmod{Gr_a}.
\]
It has a solution exactly when $a'\equiv a\pmod{m_s}$, using
\eqref{eq:ra}; then the solution for $j\bmod d_s r_a$ is unique.
Thus the entry is one precisely when $t=t'$ and
$a'\equiv a\pmod{m_s}$.

The order of $a$ modulo $m_s$ is $s$, since
$\gcd(a,m_s)=c_s$. Hence any such $a'$ has order $s$ modulo
$m_s$. If $a'\not\equiv a\pmod G$, its order modulo $d_s$
is greater than one: otherwise $d_s\mid a'$, and together with
$d_s\mid a$ and $a'\equiv a\pmod{m_s}$ this would give
$a'\equiv a\pmod G$. The two component orders are coprime,
so the order of $a'$ modulo $G$ is a proper multiple of $s$.
All off-diagonal entries therefore have strictly larger denominator.

For a bad row \eqref{eq:badrow}, all entries have the same offset
$t$, by reduction modulo $v$. Its points are
$(a+jG/p)v+t$, $0\le j<p$. For $1\le j<p$, the associated
rational number is
\[
 \frac bs+\frac jp=\frac{bp+js}{sp}.
\]
Its numerator is coprime to $s$ and nonzero modulo $p$, so its
reduced denominator is exactly $sp>s$. The $j=0$ term is the
diagonal entry, and all the points are distinct.

Finally, the $s=1$ row has only the column $(0,t)$, and the $s=2$
row has just $(0,t)$ and $(G/2,t)$. Order the matrix by offset $t$,
then by decreasing denominator among labels with $s>2$, followed
by $0$ and then $G/2$ when present. Every row has diagonal entry
one and all other entries earlier. This proves unitriangularity.
\end{proof}

\begin{corollary}\label{cor:detection}
Let $\cL\supseteq\calP$, let $x\in\W_{\cL}(Q;\F)$, and suppose
$0\le h\le\varepsilon v$ is an integer. If
$\Delta_hx(j)=x(j+h)-x(j)$ is supported on $\scrE_v$, then
$\Delta_hx=0$.
\end{corollary}
\begin{proof}
For each detecting row $C^Q_{q,u}$, both it and $C^Q_{q,u+h}$
annihilate $x$. In case~(1) of Lemma~\ref{lem:detection} this follows from the margins \eqref{e:margin};
in case~(2) it follows from $q\in\calP\subseteq\cL$.
Hence
\[
 \langle C^Q_{q,u},\Delta_hx\rangle
 =\langle C^Q_{q,u+h},x\rangle-\langle C^Q_{q,u},x\rangle=0.
\]
The supported vector is in the kernel of the unitriangular
restriction matrix from Lemma~\ref{lem:detection}, so it is zero.
\end{proof}

\subsection{Proof of Proposition~\ref{sp:prop:two}}\label{sp:sec:closure}

We continue with the field $\F$, the spaces and the sequences of \S\ref{s:finitefield} and \S\ref{s:automatic}.

\begin{proof}[Proof of Proposition~\ref{sp:prop:two}]
By Lemma \ref{lem:E_rational}, the set $E$ in \eqref{e:E_def} is a finite subset of $\Q/\Z$.
Choose a positive integer $G$ divisible by $N\prod_{p\in\calP}p$ and by the denominator of every element of $E$.
Let $H_0$ and $\varepsilon$ be as in Lemma \ref{lem:detection} for this $G$.

To apply Lemma~\ref{lem:rational_concentration}, we set $\delta = \frac{\varepsilon}{2G}$.
By Lemma~\ref{lem:rational_concentration}, for every sufficiently large $a\geq 0$ and $Q_* = Q_0 (q_1q_2)^a$, every $x\in \W_{Q_*}$ satisfies
\[x(j+P)=x(j) \qquad\text{ for all } 0\leq j<Q_* \text{ with } \dist(j/Q_*, E) \geq \delta.\]
We choose $a$ sufficiently large so that
\begin{equation}\label{e:choice_delta}
Q_* > \max\left\{\frac{1}{\delta}, \frac{GP}{\varepsilon}, P\sqrt{N}\right\}.
\end{equation}
Choose $Q\in \N$ such that $\lcm(G, Q_*)\mid Q$, i.e., there exist $r, v\in \N$ such that $Q = rQ_* = Gv$.
We assume that $v\geq 1/\varepsilon$ and $H_0\mid v$, as in Lemma \ref{lem:detection}. Such a common multiple exists, for example a sufficiently large multiple of $Q_*GH_0$.

Let $x\in \W_Q$.
We first show that $\Delta_{rP}x$ is supported on the set $\scrE_v$ of Lemma \ref{lem:detection}.
Any $i \bmod{Q}$ can be written as $i\equiv rj+t\bmod{Q}$ for $j\bmod{Q_*}$ and $0 \leq t < r$.
By Lemma \ref{lem:maps}, we have $\Res_{r,t}^{Q}x \in \W_{Q_*}$, and
\[(\Delta_P \Res_{r,t}^{Q} x)(j) = x(rj+t+rP)- x(rj+t) = (\Delta_{rP}x)(rj+t). \]
So if $(\Delta_{rP}x)(rj+t)\neq 0$, then $(\Delta_P \Res_{r,t}^{Q} x)(j)\neq 0$ and hence $\dist(j/Q_*, E)<\delta$.
Since $\frac{rj+t}{Q} = \frac{j}{Q_*} + \frac{t}{Q}$ and $\frac{t}{Q} = \frac{t}{rQ_*} < \frac{1}{Q_*} < \delta$ by \eqref{e:choice_delta}, we get
\[\dist((rj+t)/Q, E) \leq \dist(j/Q_*, E) + \delta < 2\delta.\]
Recall that $G$ is divisible by the denominator of any element in $E$.
So there exist $\frac{a_E}{G}\in E$ ($a_E\in \Z$ and not necessarily a reduced form) and $n\in \Z$ such that
\[\left|\frac{rj+t}{Q}-\frac{a_E}{G} +n \right| = \left|\frac{rj+t-a_Ev}{Q}+n\right| < 2\delta.\]
Multiplying by $Q=Gv$, the integer $t' = rj+ t-a_Ev+nQ$ satisfies $|t'| < 2\delta Q = \varepsilon v$ and $rj+t\equiv a_Ev+t'\bmod{Q}$, that is, $rj+t \in \scrE_v$.

Next, since $Q_* > \frac{GP}{\varepsilon}$,
\[rP = \frac{QP}{Q_*} = \frac{GvP}{Q_*} <  \varepsilon v.\]
So Corollary~\ref{cor:detection}, applied with $h=rP$, gives $\Delta_{rP}x =0$.

Finally, the function $x$ on $\Z$ has period $rP$ and period $Q$, hence period $g=\gcd(Q, rP)$, and $g \leq rP < \frac{rQ_*}{\sqrt{N}} = \frac{Q}{\sqrt{N}}$ by \eqref{e:choice_delta}.
By Lemma~\ref{lem:zeros}, we have $x(j) = 0$ for every integer $j$ with $-\frac{Q}{\sqrt{N}} < j < \frac{Q}{\sqrt{N}}-1$, in particular for $0\leq j\leq g-1$.
So $x=0$.
Therefore $\W_Q=0$, hence $\W_{\cL}(Q; \Q)=0$ by Lemma \ref{sp:lem:field}, and $\cL$ is vanishing.

\end{proof}

\section{Proof of Theorem~\ref{thm:main1}}\label{sec:main-proof}
The proof has four steps. First, the spanning theorem gives a reconstruction formula and a uniform dimension bound for each fixed Hecke eigenspace. Second, a compact seminorm estimate shows that $T_p-\lambda$ has finite-dimensional kernel and closed range on a suitable Banach space. Third, a second prime $q$ rules out a persistent family of eigenspaces: if the generic dimension were $n_0 > 0$, the determinant of $T_q$ on those eigenspaces would be a nonzero polynomial of degree $n_0\log q/\log p$, which is impossible for distinct primes $p$ and $q$. The possible eigenvalues are therefore discrete. Finally, the fixed growth exponent bounds the eigenvalues arising from $E_\sigma(N,k)$, and linear independence gives finiteness.

\subsection{The space \texorpdfstring{$\calH_{\vecepsilon}(N, k)$}{H\_eps(N,k)}}

Let $\vecepsilon = (\epsilon_q)_{q\in\N}$ with each $\epsilon_q \in \{\pm1\}$. 
Let $\calH_{\vecepsilon}(N, k)$ be the space of holomorphic functions $\phi(z)$ on $\H$ such that $\phi(z+1) = \phi(z)$ and $\phi(x+iy)\to 0$ as $y\to \infty$ uniformly in $x$, and for every $q\in \N$, 
\begin{equation}\label{dh:eq:FE}
\phi_{(q)}(z)= \epsilon_q\bigl(i\sqrt{\ell_N(q)}z\bigr)^{-k}
\phi_{(q)}\!\left(-\frac1{\ell_N(q)z}\right),
\quad\text{where } \phi_{(q)}(z)=\sum_{\substack{a\bmod q\\(a,q)=1}}\phi\!\left(z+\frac{a}{q}\right).
\end{equation}

Let $\phi:\H\to\C$ with $\phi(z)=\phi(z+1)$ for all $z\in\H$.
For a prime $p$, we define the normalized Hecke operator
\begin{equation}\label{e:Hecke}
(T_p \phi)(z) \coloneq p^{-\frac{k+1}{2}} \sum_{b\bmod{p}} \phi\!\left(\frac{z+b}{p}\right) + \chi_0(p) p^{\frac{k-1}{2}} \phi(pz).
\end{equation}
A standard calculation shows that the Hecke operators commute, i.e.\
$T_pT_q\phi=T_qT_p\phi$ for primes $p$, $q$. 

\begin{lemma}\label{lem:Tp_acts_H}
Let $N\in \N$, $k\in 2\N$ and $\vecepsilon=(\epsilon_q)_{q\in \N}\in\{\pm1\}^\N$.
Let $p\nmid N$ be a prime number, and suppose that $\epsilon_{pq} = \epsilon_q$ for all $q\in \N$.
Then for any $\phi\in \calH_{\vecepsilon}(N, k)$, we have $T_p\phi\in \calH_{\vecepsilon}(N, k)$. 
\end{lemma}

\begin{proof}
By definition, $(T_p \phi)(z)$ is holomorphic on $\H$ and $(T_{p}\phi)(x+iy)\to 0$ as $y\to 
\infty$ uniformly in $x$, since it inherits the property from $\phi$. 
We can also easily check that $(T_p\phi)(z+1) = (T_p\phi)(z)$.

We now check $T_p\phi$ satisfies \eqref{dh:eq:FE}. 
For each $q\in \N$, 
\begin{multline*}
(T_p \phi)_{(q)}(z) = \sum_{\substack{a\bmod{q}\\ \gcd(a, q)=1}} (T_p \phi)\!\left(z+\frac{a}{q}\right)
\\ = p^{-\frac{k+1}{2}} \sum_{\substack{a\bmod{q}\\ \gcd(a, q)=1}} \sum_{b\bmod{p}} \phi\!\left(\frac{z}{p}+\frac{a+bq}{pq}\right) 
+ p^{\frac{k-1}{2}}\sum_{\substack{a\bmod{q}\\ \gcd(a, q)=1}} \phi\!\left(pz+\frac{pa}{q}\right)
\\ = p^{-\frac{k+1}{2}} \phi_{(pq)}(z/p) 
+ \delta_{p\nmid q} \bigg\{p^{-\frac{k+1}{2}} \phi_{(q)}(z/p) 
+ p^{\frac{k-1}{2}} \phi_{(q)}(pz)\bigg\}
+ \delta_{p\mid q} \frac{\#(\Z/q\Z)^\times}{\#(\Z/(q/p)\Z)^\times} p^{\frac{k-1}{2}} \phi_{(q/p)}(pz).
\end{multline*}

Since $p\nmid N$, we have $\ell_{N}(pq) = \lcm(N, p^2q^2) = p^2 \ell_N(q)$. 
So we get
\begin{multline*}
\phi_{(pq)}(z/p)
= \epsilon_{pq} (i \sqrt{\ell_N(pq)}z/p)^{-k} \phi_{(pq)}\!\left(-\frac{1}{\ell_N(pq) z/p}\right)
\\ = \epsilon_q\bigl(i\sqrt{\ell_N(q)} z\bigr)^{-k} \phi_{(pq)}\!\left(\frac{1}{p}\left(-\frac{1}{\ell_N(q)z}\right)\right).
\end{multline*}
When $p\nmid q$, by \eqref{dh:eq:FE}, 
\begin{multline*}
p^{-\frac{k+1}{2}} \phi_{(q)}(z/p) 
+ p^{\frac{k-1}{2}} \phi_{(q)}(pz)
\\ = \epsilon_q(i\sqrt{\ell_N(q)} z)^{-k} 
\bigg\{
p^{\frac{k-1}{2}} \phi_{(q)}\!\left(-\frac{1}{\ell_N(q)z/p}\right)
+ p^{-\frac{k+1}{2}} \phi_{(q)}\!\left(-\frac{1}{\ell_N(q) pz}\right)\bigg\}
\\ = \epsilon_q(i\sqrt{\ell_N(q)} z)^{-k} 
\bigg\{
p^{\frac{k-1}{2}} \phi_{(q)}\!\left(p\left(-\frac{1}{\ell_N(q)z}\right)\right)
+ p^{-\frac{k+1}{2}} \phi_{(q)}\!\left(\frac{1}{p}\left(-\frac{1}{\ell_N(q)z}\right)\right)\bigg\}.
\end{multline*}
When $p\mid q$, since $p\nmid N$, we have $\ell_N(q) = \lcm(N, q^2p^{-2}\cdot p^2) = p^2\ell_N(q/p)$. 
Moreover $\epsilon_{q} = \epsilon_{q/p}$, so we have
\begin{multline*}
\phi_{(q/p)}(pz) = \epsilon_{q/p}\bigl(i\sqrt{\ell_N(q/p)}pz\bigr)^{-k} \phi_{(q/p)}\!\left(-\frac{1}{\ell_N(q/p)pz}\right)
\\ = \epsilon_q\bigl(i\sqrt{\ell_N(q)}z\bigr)^{-k}
\phi_{(q/p)}\!\left(p \left(-\frac{1}{\ell_N(q) z}\right)\right).
\end{multline*}
Combining the above, we obtain
\[(T_p\phi)_{(q)}(z) = \epsilon_q\bigl(i\sqrt{\ell_N(q)}z\bigr)^{-k} (T_p\phi)_{(q)}\!\left(-\frac{1}{\ell_N(q)z}\right).\]
Therefore $T_p\phi\in \calH_{\vecepsilon}(N, k)$. 
\end{proof}

\subsection{Functional equations of \texorpdfstring{$c_q$}{c\_q}-twisted \texorpdfstring{$L$}{L}-functions and \texorpdfstring{$\calH_\vecepsilon(N,k)$}{H\_eps(N,k)}}

For a sequence $(f_n)_{n \in \N}$ with $f_n = O_f(n^{\sigma_f})$, define the $c_q$-twisted $L$-series as
\begin{equation}
L_f(s, c_q) = \sum_{n=1}^\infty \frac{f_n c_q(n)}{n^s} \qquad\text{ for } \Re(s)>\sigma_f+1
\end{equation}
and the completed $c_q$-twisted $L$-function as
\begin{equation}
\Lambda_f(s, c_q) = \ell_N(q)^{\frac12(s-\frac{1}{2})} \Gamma_\C\!\left(s+\frac{k-1}{2}\right) L_f(s, c_q).
\end{equation}
Note that $L_f(s) = L_f(s, c_1)$. 

We translate the twisted functional equation into a statement about a candidate solution function $f$ reconstructed from the Dirichlet coefficients. 
For a sequence $(f_n)_{n\in\N}$ with $f_n=O_f(n^{\sigma_f})$ and $z\in\H$, define
\[
f(z) = \sum_{n=1}^\infty f_n n^{\frac{k-1}{2}} e^{2\pi inz}
\]
(converging absolutely and locally uniformly on $\H$). 
With the notation of \eqref{dh:eq:FE} applied to $\phi=f$, we have $f_{(q)}(z) = \sum_{n=1}^\infty f_n n^{\frac{k-1}{2}} c_q(n) e^{2\pi inz}$ for every $q\in\N$.

\begin{lemma}\label{lem:Euler-Ramanujan}
Given $(f_n)_{n\in \N} \in\E(N,k)$, let
$\vecepsilon=(\epsilon_q)_{q\in\N}$, where
\begin{equation}\label{e:Euler-Ramanujan-e}
\epsilon_q = \epsilon_f\prod_{\substack{p\mid(q,N)\\p^2\nmid N}} \bigl(-f_p\sqrt{p}\bigr).
\end{equation}
Then, for every $q\geq 1$, the completed $c_q$-twisted $L$-function $\Lambda_f(s,c_q)$ continues to an entire function of finite order and satisfies the functional equation
\begin{equation}\label{e:fe_hypothesis}
\Lambda_f(s,c_q) = \epsilon_q \Lambda_f(1-s,c_q).
\end{equation}
Furthermore $f\in \calH_\vecepsilon(N, k)$.
\end{lemma}

\begin{proof}
Fix $q$ and consider the twisted functional equation for $c_q$. For $p \nmid q$, the Euler factor of the twisted $L$-series at $p$ is the same as that of $L_f(s)$. For other places, we can compute the Euler factors by hand.
For each prime $p\mid q$, let $\alpha=\ord_p(q)$, and let
\[D_{f, p}(s, c_q) = (1-f_p p^{-s} + \chi_0(p)p^{-2s}) \sum_{j=0}^\infty f_{p^j} c_{p^{\alpha}}(p^j) p^{-js}.\]  Euler product factorization implies Hecke multiplicative relations for $f_n$, and \[c_q(n)=\prod_{p\mid q}c_{p^{\alpha}}\bigl(p^{\ord_p(n)}\bigr).\] 
Then
\[L_f(s, c_q) = L_f(s) \prod_{p\mid q} D_{f, p}(s, c_q)\] (see\ \cite[Lemma 4.12]{src}).

If $p\mid q$ but $p\nmid N$, we have 
\[D_{f, p}(s, c_q)
= p^{\alpha (1-s)}\big(f_{p^\alpha}  + f_{p^{\alpha-2}} p^{-1} -f_{p^{\alpha-1}} (p^{-1+s}+p^{-s})\big)
= p^{\alpha(1-2s)} D_{f, p}(1-s, c_q)\]
(by convention, $f_{p^j}=0$ when $j<0$).

Assume now that $p\mid \gcd(q, N)$. 
From the Euler product, $f_{p^j} = f_p^j$ for $p\mid N$. Thus 
\begin{multline*}
D_{f, p}(s, c_q) 
= (1-f_p p^{-s}) \bigg\{- p^{(\alpha-1)(1-s)} f_p^{\alpha-1} 
+ \varphi(p^{\alpha}) \sum_{j=\alpha}^\infty f_p^j p^{-js}\bigg\}
= p^{(\alpha-1)(1-s)}f_p^{\alpha-1}
(f_pp^{1-s}-1).
\end{multline*}
When $p^2\mid \gcd(q, N)$ we have $D_{f, p}(s, c_q)=0$, and consequently $L_f(s, c_q)=0$.
When $p\parallel N$, since $pf_p^2=1$, we have 
\[D_{f, p}(s, c_q)
= - f_p\sqrt{p} p^{(2\alpha-1)(\frac{1}{2}-s)} D_{f, p}(1-s, c_q).\]
When $p^2\mid N$ and $p\parallel q$, we have 
\[D_{f, p}(s, c_q) = f_p p^{1-s}-1 = -1.\]
Combining these together, when $\gcd(q, N)$ is squarefree, we have 
\begin{align*} 
\Lambda_f(s, c_q)
& = \Lambda_f(s) \left(\frac{\ell_N(q)}{N}\right)^{\frac{s-\frac{1}{2}}{2}} \prod_{p\mid q} D_{f, p}(s, c_q)
= \Lambda_f(s) \prod_{p\mid q}p^{(s-\frac{1}{2})\max\{0, \alpha - \frac{\ord_p(N)}{2}\}} D_{f, p}(s, c_q)
\\ &= \epsilon_f \Lambda_f(1-s) 
\cdot \prod_{\substack{p\mid \gcd(q, N)\\ p\parallel N}} \bigl(-f_p\sqrt{p}\bigr)
\cdot \prod_{p\mid q} p^{(\frac{1}{2}-s)\max\{0, \alpha-\frac{\ord_p(N)}{2}\}} D_{f, p}(1-s, c_q)
\\ &= \epsilon_q \Lambda_f(1-s, c_q).
\end{align*}

The displayed factorization expresses $\Lambda_f(s,c_q)$ as $\Lambda_f(s)$ times an entire function of exponential type, so it is entire of finite order. Absolute convergence and Stirling's formula give boundedness on a sufficiently far right vertical line, and the functional equation gives boundedness on the reflected line. Since the completed function is entire of finite order, Phragm\'en--Lindel\"of gives boundedness on every vertical strip. Hecke's theorem therefore applies.
By Hecke's argument \cite[Theorem~4.3.5]{Miyake}, the functional equation \eqref{e:fe_hypothesis} gives \eqref{dh:eq:FE} for $f$.
Its Fourier series makes $f$ holomorphic and periodic on $\H$, with $f(x+iy)\to0$ uniformly in $x$ as $y\to\infty$. Hence $f\in\calH_{\vecepsilon}(N,k)$.
\end{proof}

\begin{lemma}\label{lem:Hecke-Euler}
For $(f_n)\in\E(N,k)$ and every prime $p$, we have $T_pf=f_pf$.
\end{lemma}
\begin{proof}
Set $f_x=0$ when $x\notin\N$. The Euler product gives multiplicativity and
\[
 f_{p^{j+1}}=f_pf_{p^j}-\chi_0(p)f_{p^{j-1}}\qquad(j\geq0).
\]
Consequently $f_{pn}+\chi_0(p)f_{n/p}=f_pf_n$ for every $n\geq1$.
Applying \eqref{e:Hecke} term by term to the Fourier series therefore gives
\[
 (T_pf)(z)=\sum_{n\geq1}\bigl(f_{pn}+\chi_0(p)f_{n/p}\bigr)
 n^{\frac{k-1}{2}}e^{2\pi inz}=f_pf(z).
\]
\end{proof}

\begin{lemma}\label{lem:linindep}
The set $\E(N,k)$ is linearly independent over $\C$.
\end{lemma}
\begin{proof}
Distinct solutions have distinct tuples $(f_p)_p$, since the Euler product determines every coefficient from the coefficients at primes.
Suppose there is a nontrivial relation $\sum_{i=1}^r c_i f^{(i)}=0$ among distinct solutions, with $r$ minimal and every $c_i\neq0$.
Since $f^{(i)}_1=1$, we have $r\geq2$.
Choose a prime $p$ with $f^{(1)}_p\neq f^{(2)}_p$ and apply $T_p-f^{(1)}_p$ to the relation.
By Lemma~\ref{lem:Hecke-Euler}, the first term disappears and the second remains nonzero, giving a nontrivial relation with fewer than $r$ terms, a contradiction.
\end{proof}

\begin{theorem}\label{thm:main-general}
Let $N,k\in\N$ with $k$ even, and let $\vecepsilon = (\epsilon_q)_{q\in \N}\in\{\pm1\}^\N$. 
Let $\calP$ be a finite set of primes which is vanishing for $N$, and define
\begin{equation}\label{e:calH_T-eigen}
\calH_{\vecepsilon}(N,k;\calP,\veclambda)
\coloneq \left\{\phi\in \calH_{\vecepsilon}(N,k): T_p \phi = \lambda_{p} \phi \; \forall p\in \calP\right\}
\quad\text{for }\veclambda=(\lambda_p)_{p\in\calP}\in\C^\calP.
\end{equation}
For any fixed tuple $(M_p)_{p\in\calP}$, we have $\dim\calH_{\vecepsilon}(N,k;\calP,\veclambda) = O_{N,\calP,(M_p)}(k)$, uniformly in $\vecepsilon$ and $\veclambda$ satisfying $|\lambda_p|\le M_p$ for all $p\in\calP$. 
\end{theorem}

Assume that some set $\calP$ of primes dividing $N$ is vanishing for $N$. 
When $(f_n)_{n\in \N}\in \E(N, k)$, by Lemma \ref{lem:Euler-Ramanujan} and Lemma \ref{lem:Hecke-Euler}, we have $f\in \calH_{\vecepsilon}(N,k;\calP,\veclambda)$, with $\epsilon_q$ as given by \eqref{e:Euler-Ramanujan-e} and $\lambda_p = f_p$ with $\sqrt{p} f_p\in \{0,\pm1\}$ for all $p\in\calP$.
So both $\vecepsilon$ and $\veclambda$ are determined by $(f_n)_{n\in \N}$ 
and the elements in $\E(N, k)$ lie in at most $2^{\omega(N)+1}$ possible spaces $\calH_{\vecepsilon}(N,k;\calP,\veclambda)$. 
By Lemma \ref{lem:linindep} and Theorem \ref{thm:main-general}, we have 
\[\#\E(N, k) \leq \sum_{\vecepsilon} \dim\calH_{\vecepsilon}(N,k;\calP,\veclambda)
= O_N(k).\]
This proves the following corollary. 
\begin{corollary}\label{cor:main-dividingN}
Whenever some set of primes dividing $N$ is vanishing for $N$, 
\[\#\E(N, k) = O_N(k).\]
\end{corollary}
By part~(2) of Theorem~\ref{thm:spanning}, the set of primes dividing $N$ is vanishing for $N$ when $\rad(N)>\sqrt N$, so Corollary~\ref{cor:main-dividingN} gives the finiteness assertion of Theorem~\ref{thm:main-asymptotic}.

\subsection{Proof of Theorem \ref{thm:main-general}}\label{ss:proof_main-general}

Let $N$ be a positive integer. 
We fix $\calP$, a finite set of primes which is vanishing for $N$, and nonnegative bounds $M_p$ for $p\in\calP$. We allow any $\veclambda=(\lambda_p)_{p\in\calP}$ with $|\lambda_p|\le M_p$. 
Choose a positive integer $Q$ such that $\calR_{N}(Q) \cup \bigl\{C_{p, u}^Q:\, p\in \calP, \, u\bmod{Q/p}\bigr\}$ spans $\Q^{\Z/Q\Z}$. We regard these data as fixed for the remainder of the proof.

For a function $f\in \calH_{\vecepsilon}(N,k;\calP,\veclambda)$, we write the Fourier expansion as 
\[f(z) = \sum_{n\geq 1} f_n n^{\frac{k-1}{2}} e^{2\pi inz}. \]
To get an upper bound for the dimension of $\calH_{\vecepsilon}(N,k;\calP,\veclambda)$, it suffices to show that there exists a $B=O_{N,\calP,(M_p)}(k)$ such that $f=0$ whenever $f_n=0$ for $n\leq B$. 

We fix choices of coefficients $A_j(q, u)$ and $B_j(p, m)$ such that 
\begin{equation}\label{e:certificate}
e_j = \sum_{q, u} A_j (q, u) R_{q, u}^Q + \sum_{p, m} B_j(p, m) C_{p, m}^Q \quad \text{for each }j\bmod Q,
\end{equation}
where the first sum runs over $\calR_{N}(Q)$ and the second sum runs over $p\in \calP$ and $m\bmod{Q/p}$. 

Define 
\begin{equation}\label{e:cost}
C = \max\left(\{1\}\cup\bigg\{\sum_{q, u} |A_{j}(q, u)|\varphi(q) + \sum_{p, m} |B_{j}(p, m)| (\sqrt{p}M_p + \chi_0(p)):j\bmod Q\bigg\}\right).
\end{equation}
By the definition of legal rows \eqref{e:scrR_NQ_def}, we can choose $0< \theta< 1$ such that 
\begin{equation}\label{e:buffer}
\frac{\sqrt{\ell_N(q)}}{Q} \max\{|u|, |u+1|\} \leq \theta \qquad\text{ whenever } A_{j}(q, u)\neq 0. 
\end{equation}
Let
\begin{equation}\label{e:constants}
y_0 = \sqrt{\frac{1-\theta^2}{2\ell_N(Q)}}, \qquad c(\theta) = \frac{2}{1+\theta^2}>1.
\end{equation}
Fix $\nu\geq k/2$.
Lemmas~\ref{lem:reconstruction} and~\ref{lem:envelope} hold for any such choices; each later application specifies $\calP$, the bounds $M_p$ and $\nu$, and makes the choices above for them.
For $f:\H\to\C$, we define
\begin{equation}\label{e:sup_f_def}
\|f\|_{\infty, \nu} \coloneq \sup_{z\in \H} \Im(z)^\nu |f(z)|,
\end{equation}
and
\begin{equation}\label{e:calH_nu_def}
\calH_{\vecepsilon,\nu}(N,k) =\{f\in\calH_{\vecepsilon}(N,k):\|f\|_{\infty,\nu}<\infty\}.
\end{equation}

\begin{lemma}\label{lem:reconstruction}
For $x\in [j/Q, (j+1)/Q]$ and $y>0$, 
let $\beta_j(q, u) = x+iy+\frac{u-j}{Q}$. 
For $f\in \calH_{\vecepsilon}(N, k)$, 
\begin{multline}\label{e:reconstruction}
f(z) = \sum_{q, u} A_{j}(q, u) \epsilon_q (i\sqrt{\ell_N(q)}\beta_j(q, u))^{-k} 
\sum_{\substack{a\bmod{q}\\\gcd(a, q)=1}} f\!\left(-\frac{1}{\ell_N(q)\beta_j(q, u)}+\frac{a}{q}\right)
\\ + \sum_{p, m} B_{j}(p, m) \bigg\{p^{\frac{k+1}{2}} (T_p f)(p \beta_j(p, m)) - \chi_0(p) p^{k} f(p^2\beta_j(p, m))\bigg\}.
\end{multline}
For $0 < y \leq y_0$ and $A_j(q,u)\neq0$, we have 
\[\Im\left(-\frac{1}{\ell_N(q)\beta_j(q, u)}+\frac{a}{q}\right) 
= \frac{y}{\ell_N(q)|\beta_j(q, u)|^2}
\geq y c(\theta).\]
\end{lemma}

\begin{proof} 
We define $v\in \C^{\Z/Q\Z}$ via 
\[v(n) = \left<v, e_n\right> = f\!\left(z+\frac{n-j}{Q}\right). \]
Applying \eqref{e:certificate}, we have 
\begin{align*}
f(z) & = v(j) = \left<v, e_j\right> = \sum_{q, u} A_{j}(q, u) \sum_{\substack{a\bmod{q}\\ \gcd(a, q)=1}} v(u+aQ/q) 
+ \sum_{p, m}B_j(p,m) \sum_{b\bmod{p}} v(m+bQ/p)
\\ & = \sum_{q, u} A_{j}(q, u) \sum_{\substack{a\bmod{q}\\ \gcd(a, q)=1}} f\!\left(z+\frac{u-j}{Q}+\frac{a}{q}\right) 
+ \sum_{p, m} B_{j}(p, m)
\sum_{b\bmod{p}} f\!\left(z+\frac{m-j}{Q}+\frac{b}{p}\right).
\end{align*}
Recall that $\sum_{\substack{a\bmod{q}\\ \gcd(a, q)=1}} f(z+a/q) = f_{(q)}(z)$ and 
$\sum_{b\bmod{p}} f(z+b/p) = p^{\frac{k+1}{2}} (T_p f)(p z) - p^{k} \chi_0(p) f(p^2 z)$. 
We have 
\begin{multline*}
f(z) = \sum_{q, u} A_{j}(q, u) f_{(q)}\!\left(z+\frac{u-j}{Q}\right) 
\\ + \sum_{p, m} B_{j}(p, m) \bigg\{p^{\frac{k+1}{2}} (T_p f)\!\left(p \left(z + \frac{m-j}{Q}\right)\right) 
- \chi_0(p) p^{k} f\!\left(p^2 \left(z + \frac{m-j}{Q}\right)\right) \bigg\}. 
\end{multline*}

Write $z = x+iy$, $y>0$ and $x\in [j/Q, (j+1)/Q]$. 
Recall that $\beta_j(q, u) = x+iy+\frac{u-j}{Q}$. 
By applying \eqref{dh:eq:FE}, for $f_{(q)}$, we get \eqref{e:reconstruction}.
By our choice of $x$, we have 
\[x+\frac{u-j}{Q} \in \left[\frac{u}{Q}, \frac{u+1}{Q}\right].\]
So by using the constant given in \eqref{e:buffer}, 
\[\left|x+\frac{u-j}{Q}\right| \leq \frac{\theta}{\sqrt{\ell_N(q)}} \qquad\text{ whenever } A_{j}(q, u)\neq 0, \]
and we have 
\[\ell_N(q) |\beta_j(q, u)|^2 \leq \theta^2 + \ell_N(q) y_0^2 \leq \theta^2 + \ell_N(Q) y_0^2=\frac{1+\theta^2}{2} = \frac{1}{c(\theta)}  \]
for $0 < y \leq y_0$, since $\ell_N(q)$ divides $\ell_N(Q)$.
We have 
\[\Im\left(-\frac{1}{\ell_N(q)\beta_j(q, u)}+\frac{a}{q}\right) 
= \frac{y}{\ell_N(q)|\beta_j(q, u)|^2}
\geq y c(\theta).\]
\end{proof}

The next lemma bounds $\|f\|_{\infty,\nu}$ by the values of $f$ at height at least $y_0$. 
Its second part requires $T_pf=\lambda_pf$ only up to the error $(T_p-\lambda_p)f$; Lemma~\ref{lem:closed-range} uses this to prove that $T_p-\lambda$ has closed range.

\begin{lemma}\label{lem:envelope}
With $Q$, the coefficients $A_j(q,u)$ and $B_j(p,m)$, and $\theta$, $y_0$ and $c(\theta)$ as in \eqref{e:certificate}--\eqref{e:constants}, define
\[\rho_\nu =\max_{j\bmod{Q}}\left\{c(\theta)^{\frac{k}{2}-\nu} \sum_{q, u} |A_{j}(q, u)| \varphi(q) + \sum_{p\in\calP} \sum_{m} |B_j(p, m)| \bigl(M_pp^{\frac{k+1}{2}-\nu}+\chi_0(p)p^{k-2\nu}\bigr)\right\}, \]
and for each $p\in \calP$, 
\[D_{p,\nu} = p^{\frac{k+1}{2}-\nu} \max_{j\bmod{Q}} \sum_{m} |B_j(p, m)|.\]
Then $\rho_\nu\le Cc(\theta)^{\frac{k}{2}-\nu}$. 
If $\rho_\nu\le1$ and
$f\in\calH_{\vecepsilon}(N,k;\calP,\veclambda)$, where
$|\lambda_p|\le M_p$, then
\begin{equation}\label{e:envelope}
 \|f\|_{\infty,\nu}=\sup_{x\in\R,\ y\ge y_0}y^\nu|f(x+iy)|<\infty.
\end{equation}
If $\rho_\nu<1$, then every $f\in\calH_{\vecepsilon,\nu}(N,k)$
and every such tuple $\veclambda$ satisfy
\begin{equation}\label{e:residual-envelope}
\|f\|_{\infty,\nu}
\leq 
\max\left\{\sup_{x\in\R,\ y\geq y_0} y^\nu|f(x+iy)|, \frac{1}{1-\rho_\nu} \sum_{p\in\calP}D_{p,\nu} \|(T_p-\lambda_p)f\|_{\infty,\nu}\right\}.
\end{equation}
Here $\|(T_p-\lambda_p)f\|_{\infty,\nu}$ is defined by \eqref{e:sup_f_def} even when $T_pf\notin\calH_{\vecepsilon}(N,k)$.
\end{lemma}

\begin{proof}
Write $G(e^{2\pi iz})=f(z)$. 
Periodicity and uniform vanishing at
infinity imply that $G$ is holomorphic in the unit disk with $G(0)=0$.
The maximum principle for $G(w)/w$ gives, for $y \geq a > 0$,
\begin{equation}\label{e:cusp-tail}
\sup_{x\in \R} |f(x+iy)|\le e^{-2\pi(y-a)}\sup_{x\in \R} |f(x+ia)|.
\end{equation}
In particular 
\begin{equation}\label{e:sup_up_finite}
\sup_{x\in\R,\ y\ge y_0}y^\nu|f(x+iy)|<\infty, 
\end{equation}
without a growth assumption near
the real axis.

For any primes $p$, note that $(T_p f)(z)\to 0$ as $\Im(z)\to \infty$ uniformly in $\Re(z)$. 
By \eqref{e:reconstruction}, for $\Re(z)\in [j/Q, (j+1)/Q]$ and $\Im(z) = y\in (0, y_0]$, we have 
\begin{multline*}
y^{\nu} |f(z)| \leq \sum_{q, u} |A_{j}(q, u)| \varphi(q)
(y{y'}^{-1})^{\nu-\frac{k}{2}} \max_{x\in\R} {y'}^{\nu} |f(x+iy')| 
\\ + \sum_{p, m} |B_{j}(p, m)| \bigg(p^{\frac{k+1}{2}-\nu} \max_{x\in \R} (py)^{\nu}|(T_p-\lambda_p)f(x+ipy)| + p^{\frac{k+1}{2}-\nu} |\lambda_p| \max_{x\in\R} (py)^{\nu} |f(x+ipy)| 
\\ + \chi_0(p) p^{k-2\nu} \max_{x\in \R}(p^2y)^{\nu} |f(x+ip^2y)|\bigg). 
\end{multline*}
Here, in the sum over $q, u$, $y' = \Im\left(-\frac{1}{\ell_N(q)\beta_j(q, u)}\right) = \frac{y}{\ell_N(q)|\beta_j(q, u)|^2} \geq yc(\theta)$. 
Since $\nu-\frac{k}{2}\geq0$, we have $(y{y'}^{-1})^{\nu-\frac{k}{2}} \leq c(\theta)^{\frac{k}{2}-\nu}$.
Also note that $p \geq 2 > c(\theta) >1$.
Let 
\[a_j=\sum_{q,u}|A_j(q,u)|\varphi(q), 
\qquad b_{j,p}=\sum_m|B_j(p,m)| \qquad \text{ for } j\bmod{Q}.\]
We have 
\begin{multline*}
y^{\nu} |f(z)| \leq \bigg\{a_j c(\theta)^{\frac{k}{2}-\nu} + \sum_{p\in \calP} b_{j, p} (p^{\frac{k+1}{2}-\nu} |\lambda_p| + \chi_0(p)p^{k-2\nu})\bigg\}\max_{x\in \R, y'\geq yc(\theta)} {y'}^{\nu} |f(x+iy')|
\\ + \sum_{p\in \calP} b_{j, p} 
p^{\frac{k+1}{2}-\nu} \max_{x\in \R, y'\geq yc(\theta)} {y'}^{\nu} |((T_p-\lambda_p) f) (x+iy')|.
\end{multline*}
Taking the maximum over $j\bmod{Q}$ and supremum for ${y'}^{\nu}|f(x+iy')|$ for $y'\geq c(\theta)y$, we get
\begin{multline}\label{e:one-step}
y^{\nu} |f(x+iy)| 
\leq \rho_{\nu} \sup_{x\in \R, y'\geq yc(\theta)} {y'}^{\nu} |f(x+iy')|
\\ + \sum_{p\in \calP} D_{p, \nu} \sup_{x\in \R, y'\geq yc(\theta)} {y'}^{\nu} |((T_p-\lambda_p) f) (x+iy')|.
\end{multline}
Note that the definition of $C$ gives $\rho_\nu \leq Cc(\theta)^{\frac{k}{2}-\nu}$. 
When $f\in \calH_{\vecepsilon}(N,k;\calP,\veclambda)$, the terms involving $(T_p-\lambda_p)f$ vanish.
Assume that $\rho_\nu\leq 1$.
Expand each value of $f$ at height below $y_0$ on the right of \eqref{e:reconstruction} by \eqref{e:reconstruction} again, and repeat.
Each step increases the height by at least $c(\theta)>1$, so for each fixed starting height this procedure stops after finitely many steps.
At every step the sum of the absolute weighted coefficients is at most $\rho_\nu\leq 1$.
Thus the final sum of the absolute weighted coefficients is at most one, and all remaining values have height at least $y_0$.
It follows that
\[y^\nu |f(x+iy)|\leq \sup_{x\in\R,\ y'\geq y_0}{y'}^\nu|f(x+iy')|.\]
Together with \eqref{e:sup_up_finite}, this proves \eqref{e:envelope}.

Now assume that $f\in\calH_{\vecepsilon,\nu}(N,k)$ and $\rho_\nu<1$.
Taking the supremum in \eqref{e:one-step} for $0<y\leq y_0$, and treating $y\geq y_0$ separately, gives
\[
\|f\|_{\infty,\nu}\leq
\max\left\{
\sup_{x\in\R,\ y\geq y_0}y^\nu|f(x+iy)|,
\rho_\nu\|f\|_{\infty,\nu}
+\sum_{p\in\calP}D_{p,\nu}\|(T_p-\lambda_p)f\|_{\infty,\nu}
\right\}.
\]
If the norm is larger than the first term, then
\[
(1-\rho_\nu)\|f\|_{\infty,\nu}
\leq\sum_{p\in\calP}D_{p,\nu}\|(T_p-\lambda_p)f\|_{\infty,\nu}.
\]
This proves \eqref{e:residual-envelope}.

\end{proof}

\begin{proof}[Proof of Theorem \ref{thm:main-general}]

Set $\nu=\frac{k}{2}+\frac{\log C}{\log c(\theta)}$ and
\begin{equation}\label{e:constants2}
B=\left\lfloor\frac{\nu\log2}{\pi y_0}\right\rfloor.
\end{equation}
Then $\rho_\nu\le Cc(\theta)^{\frac{k}{2}-\nu}=1$, so Lemma~\ref{lem:envelope} applies.
For any $f\in \calH_{\vecepsilon}(N,k;\calP,\veclambda)$, let $G(e^{2\pi iz}) = f(z)$. 
As in the proof of Lemma \ref{lem:envelope}, $G$ is holomorphic on the unit disk. 
We further assume that $f_1=\cdots=f_B=0$. 
Then $G(w)$ has a zero of order at least $B+1$ at $w=0$. 
Applying the maximum principle to $\frac{G(w)}{w^{B+1}}$ on $|w|= e^{-\pi y_0}$ gives 
\[\sup_{x\in \R} |f(x+iy)| \leq e^{-2\pi(B+1)(y-\frac{y_0}{2})} \sup_{x\in \R} |f(x+iy_0/2)|,  \qquad y\geq \frac{y_0}{2}.\]
By Lemma \ref{lem:envelope}, for any $y \geq y_0$, we have, 
\[y^\nu |f(x+iy)| \leq e^{-2\pi(B+1)(y-\frac{y_0}{2})} \left(\frac{2 y}{y_0}\right)^\nu\|f\|_{\infty, \nu}.\]
Let $I(y) = e^{-2\pi(B+1)(y-\frac{y_0}{2})} (2yy_0^{-1})^{\nu}$. 
Since 
\[\frac{d}{dy} \log(I(y)) = \frac{\nu}{y} -2\pi(B+1) \leq \frac{\nu}{y_0}(1 - 2\log 2) < 0, \]
$I(y)$ is decreasing. 
By taking $y=y_0$, we have 
\[y^{\nu} |f(x+iy)| \leq e^{-\pi(B+1) y_0} 2^\nu \|f\|_{\infty, \nu}.\]
Recalling \eqref{e:constants2}, we have $\pi(B+1) y_0 -\nu\log 2 > 0$. 
Taking the supremum and using \eqref{e:envelope} gives $\|f\|_{\infty,\nu}\le 2^\nu e^{-\pi(B+1)y_0}\|f\|_{\infty,\nu}$; the multiplier is strictly less than $1$, so $f=0$.
Therefore $\dim\calH_{\vecepsilon}(N,k;\calP,\veclambda) \leq B = O_{N,\calP,(M_p)}(k)$.
\end{proof}

\subsection{Finiteness of Hecke eigenvalues}\label{ss:proof-finite-HE}

To complete the proof of Theorem \ref{thm:main1}, it remains to show that, for a fixed prime $p>\sqrt{N}$ not dividing $N$, only finitely many eigenvalues of $T_p$ arise from $\E_\sigma(N,k)$.

We now restrict to a sign system $\vecepsilon= (\epsilon_{q})_{q\in \N}$ satisfying
\begin{equation}\label{e:good-signs}
\epsilon_{pq}=\epsilon_q \qquad 
\text{ for any prime } p\nmid N, \text{ and any } q\in \N.
\end{equation}
Any $\vecepsilon$ arising from \eqref{e:Euler-Ramanujan-e} has this
property. 
Throughout the remainder of the proof, $N$, $k$ and $\vecepsilon$ are fixed, and for $\nu>0$ we write $X_\nu\coloneq\calH_{\vecepsilon,\nu}(N,k)$.

For every $\nu>0$, the normed space $X_\nu$ is Banach.
Indeed, any Cauchy sequence $(\phi_n)_{n\in\N}$ in this space, with respect to $\|\cdot\|_{\infty,\nu}$, converges locally uniformly on $\H$.
Let $\phi$ be its limit.
Then $\phi$ is a holomorphic periodic function satisfying all of \eqref{dh:eq:FE}.
Since the sequence is norm-bounded, taking pointwise limits gives
\[
\|\phi\|_{\infty,\nu}\leq\sup_n\|\phi_n\|_{\infty,\nu}<\infty.
\]
Hence $|\phi(x+iy)|\leq y^{-\nu}\|\phi\|_{\infty,\nu}$, so $\phi(x+iy)$ tends to zero uniformly in $x$ as $y\to\infty$.
Thus $\phi\in X_\nu$.
Taking pointwise limits in the differences also gives
\[
\|\phi_n-\phi\|_{\infty,\nu}
\leq\sup_{m\geq n}\|\phi_n-\phi_m\|_{\infty,\nu}\longrightarrow0.
\]
This proves completeness.

For any prime $p\nmid N$, by Lemma~\ref{lem:Tp_acts_H}, we see that the Hecke operator $T_p$ is a bounded operator on $X_\nu$, with
\begin{equation}\label{e:Hecke-norm}
\|T_p\|_{X_\nu} := \sup_{f\in X_\nu, f\neq 0} \frac{\|T_p f\|_{\infty, \nu}}{\|f\|_{\infty, \nu}} \le p^{\nu-\frac{k-1}{2}} + p^{\frac{k-1}{2}-\nu}.
\end{equation}
In fact the first term bounds the sum over $b\bmod{p}$ in \eqref{e:Hecke} and the second bounds the remaining term in \eqref{e:Hecke}. 
The evaluations in the sum have height $\Im(z)/p$, and the remaining evaluation has height $p\Im(z)$. 

With $y_0$ as in \eqref{e:constants}, every norm-bounded sequence in $X_\nu$ has a subsequence $(\phi_j)$ such that
\begin{equation}\label{e:semi-norm-cauchy}
\sup_{x\in\R,\ y\geq y_0}y^\nu\bigl|(\phi_{j_1}-\phi_{j_2})(x+iy)\bigr| \to 0
\qquad\text{ as } j_1, j_2\to \infty.
\end{equation}
Indeed, the norm bound implies local uniform boundedness on $\H$, so Montel's theorem gives locally uniform convergence along a subsequence. 
By periodicity this controls every strip
$0\le x\le1$, $y_0\le y\le Y$. 
Also \eqref{e:cusp-tail} with
$a=y_0/2$ bounds the remaining tails, uniformly along the sequence,
by a fixed multiple of $y^\nu e^{-2\pi(y-y_0/2)}$, which tends to
zero. 
Combining the strip and tail estimates proves \eqref{e:semi-norm-cauchy}.

\begin{lemma}\label{lem:closed-range}
Let $p>\sqrt N$ be a prime not dividing $N$, and let $M<\infty$.
There is $\nu\ge k/2$, depending on $N$, $p$, $M$, $k$, such that for every
$|\lambda|\le M$, the operator $T_p-\lambda$ on $X_\nu$ has finite-dimensional kernel and closed
range. 
Moreover,
\begin{equation}\label{e:full-kernel}
\ker(T_p-\lambda:X_\nu \to X_\nu)
=\{f\in\calH_{\vecepsilon}(N,k):T_pf=\lambda f\}.
\end{equation}
\end{lemma}

\begin{proof}
We apply Theorem~\ref{thm:spanning} with $\calP=\{p\}$. 
Make the choices of \S\ref{ss:proof_main-general} with $\calP=\{p\}$ and $M_p=M$: fix $Q$ and the coefficients in \eqref{e:certificate}, and let $C$, $\theta$, $y_0$ and $c(\theta)$ be as in \eqref{e:cost}--\eqref{e:constants}. 
Choose $\nu>\frac{k}{2}+\frac{\log C}{\log c(\theta)}$.
Then Lemma~\ref{lem:envelope} gives $\rho_\nu\le Cc(\theta)^{\frac{k}{2}-\nu}<1$.
Then, by \eqref{e:envelope}, for any $f\in \calH_{\vecepsilon}(N,k)$ satisfying $T_pf=\lambda f$, we have $\|f\|_{\infty, \nu} = \sup_{x\in \R, y\geq y_0} y^{\nu} |f(x+iy)| < \infty$, which proves 
the equality \eqref{e:full-kernel}.
By Lemma \ref{lem:envelope}, for any $f\in X_\nu$, we have 
\begin{equation}\label{e:compact-estimate}
\|f\|_{\infty, \nu} \le\max\left\{\sup_{x\in\R,\ y\ge y_0}y^\nu |f(x+iy)|, \frac{D_{p,\nu}}{1-\rho_\nu} \|(T_p-\lambda)f\|_{\infty, \nu} \right\}. 
\end{equation}

For $f\in\ker(T_p-\lambda)$, \eqref{e:compact-estimate} gives $\|f\|_{\infty,\nu}=\sup_{x\in\R,\ y\ge y_0}y^\nu|f(x+iy)|$.
Every sequence in the closed unit ball of $\ker(T_p-\lambda)$ therefore has a Cauchy subsequence in the norm $\|\cdot\|_{\infty, \nu}$ by \eqref{e:semi-norm-cauchy}. 
By completeness and continuity of $T_p-\lambda$, this subsequence converges to a point in the same closed unit ball.
Thus the closed unit ball is compact, and $\ker(T_p-\lambda)$ is finite dimensional.
It therefore has a closed complement $Y$ in $X_\nu$.
If $(T_p-\lambda)|_Y$ were not bounded below, there would be $g_j\in Y$ with $\|g_j\|_{\infty,\nu}=1$ and $(T_p-\lambda)g_j\to0$ in norm.
By \eqref{e:semi-norm-cauchy}, a subsequence is Cauchy for the seminorm $\phi\mapsto\sup_{x\in\R,\ y\ge y_0}y^\nu|\phi(x+iy)|$.
Applying \eqref{e:compact-estimate} to differences shows that this subsequence is also Cauchy in $\|\cdot\|_{\infty,\nu}$.
Its limit belongs to $Y\cap\ker(T_p-\lambda)=\{0\}$ but has norm one, a contradiction.
Thus there is a constant $\kappa>0$ such that
\[\|(T_p-\lambda)g\|_{\infty,\nu}\geq \kappa\|g\|_{\infty,\nu}\qquad\text{ for any } g\in Y.\]
If a sequence $((T_p-\lambda)g_j)$ with $g_j\in Y$ converges in $X_\nu$, this inequality applied to differences shows that $(g_j)$ is Cauchy.
Since $Y$ is closed in a Banach space, $g_j\to g\in Y$, and continuity gives $(T_p-\lambda)g_j\to(T_p-\lambda)g$.
Hence $(T_p-\lambda)(Y)$ is closed.
Every element of $X_\nu$ is the sum of an element of the kernel and an element of $Y$, so
\[
(T_p-\lambda)\bigl(X_\nu\bigr)=(T_p-\lambda)(Y).
\]
This proves that $T_p-\lambda$ has closed range.
\end{proof}

For $x\geq 0$, define 
\[\log_p^+x = \frac{\log(\max\{1, x\})}{\log p}.\]

\begin{lemma}\label{lem:two-prime}
Let $p$ and $q$ be prime numbers with $(N,pq)=1$ and $p,q>\sqrt N$. 
Suppose $f\in\calH_\vecepsilon(N,k)$ is a common eigenfunction of $T_p$ and $T_q$, with eigenvalues $\lambda_p$ and $\lambda_q$.
Then
\[\log_q^+|\lambda_q|\le\log_p^+|\lambda_p|+O_{N,p}(1).\]
In particular, there are positive constants $c_{p,q},C_{p,q}$ such that
\begin{equation}\label{e:power-comparison}
 c_{p,q}|\lambda_p|^{\log q/\log p}
 \le |\lambda_q|\le
 C_{p,q}|\lambda_p|^{\log q/\log p}
 \quad\text{whenever $|\lambda_p|$ is sufficiently large}.
\end{equation}
\end{lemma}

\begin{proof}
Make the choices of \S\ref{ss:proof_main-general} with $\calP=\{p\}$, and let $\rho_\nu$ be as in Lemma~\ref{lem:envelope}. Choose $A_{N,p}\geq0$, depending only on $N$ and $p$, and set $\nu_0=k/2+A_{N,p}$ so that $\rho_{\nu_0}<1$ when $M_p=1$; this is possible because, for $M_p=1$, $\rho_\nu$ depends on $\nu$ only through $\nu-\frac{k}{2}$ and tends to $0$ as $\nu\to\infty$.
Set $M_p=|\lambda_p|$ and $\nu=\nu_0+\log_p^+|\lambda_p|$.
Since $c(\theta)^{\frac{k}{2}-\nu}\le c(\theta)^{\frac{k}{2}-\nu_0}$, $|\lambda_p|\,p^{\frac{k+1}{2}-\nu}\le p^{\frac{k+1}{2}-\nu_0}$ and $p^{k-2\nu}\le p^{k-2\nu_0}$, each term of $\rho_\nu$ is at most the corresponding term of $\rho_{\nu_0}$ with $M_p=1$, so $\rho_\nu<1$.
Applying Lemma~\ref{lem:envelope} with $\calP=\{p\}$ and $\veclambda=(\lambda_p)$ gives
$f\in X_\nu$, without any assumed growth bound.
By \eqref{e:Hecke-norm},
\[|\lambda_q|\le q^{\nu-\frac{k-1}{2}}+q^{\frac{k-1}{2}-\nu}
\le2q^{\nu-\frac{k-1}{2}}.\]
Since $q\ge2$, we have $\log_q2\le1$, and hence
\[\log_q^+|\lambda_q|\le\log_p^+|\lambda_p|+A_{N,p}+\frac32
=\log_p^+|\lambda_p|+O_{N,p}(1).\]
Interchanging $p$ and $q$ gives
\[\log_q^+|\lambda_q|=\log_p^+|\lambda_p|+O_{N,p,q}(1).\]
For sufficiently large $|\lambda_p|$, this forces $|\lambda_q|>1$, and
exponentiation gives \eqref{e:power-comparison}.
\end{proof}

Consider a prime $p>\sqrt N$ not dividing $N$. 
By Theorem \ref{thm:spanning}, $\{p\}$ is vanishing for $N$. 
For $\lambda\in \C$, let
\[K_p(\lambda) = \calH_{\vecepsilon}(N,k;\{p\},(\lambda)) = \{f\in \calH_{\vecepsilon}(N, k):\, T_p f = \lambda f\}.\]

The proof of the next proposition uses the following notions and results of M\"uller \cite{Muller2000}, which we state for the family $T_p-\lambda$ of bounded operators on the Banach space $X_\nu$, depending analytically on $\lambda$. The family is \emph{regular} at $\lambda_0$ if $T_p-\lambda_0$ has closed range and the kernels converge in the gap metric, that is,
\[\sup_{\substack{x\in\ker(T_p-\lambda_0)\\ \|x\|_{\infty,\nu}\le1}}\dist\bigl(x,\ker(T_p-\lambda)\bigr)+\sup_{\substack{x\in\ker(T_p-\lambda)\\ \|x\|_{\infty,\nu}\le1}}\dist\bigl(x,\ker(T_p-\lambda_0)\bigr)\to0\qquad\text{as }\lambda\to\lambda_0\]
(see \cite[Theorem~1.1]{Muller2000} and the definition following it). If the family is regular at $\lambda_0$, then every $x\in\ker(T_p-\lambda_0)$ extends to a holomorphic function $x(\lambda)$ near $\lambda_0$ with $x(\lambda_0)=x$ and $(T_p-\lambda)x(\lambda)=0$ \cite[Theorem~1.2]{Muller2000}. If $T_p-\lambda_0$ has closed range and finite-dimensional kernel, then the family is regular on a punctured neighborhood of $\lambda_0$, and $\dim\ker(T_p-\lambda)$ is constant there and at most $\dim\ker(T_p-\lambda_0)$ \cite[Theorems~1.6 and~1.8]{Muller2000}.

\begin{proposition}\label{prop:point-discrete}
For every sign system satisfying \eqref{e:good-signs} and every
prime $p>\sqrt N$ not dividing $N$, the set
\[
 \{\lambda\in\C:K_p(\lambda)\ne0\}
\]
has no finite accumulation point. In particular, only finitely
many such eigenvalues lie in any bounded subset of $\C$.
\end{proposition}

\begin{proof}
Fix a prime $p>\sqrt N$ not dividing $N$, and choose a distinct prime $q>\sqrt N$, $q\nmid N$. 
For any $\lambda\in \C$, 
the space $K_p(\lambda)$ is finite dimensional by Theorem~\ref{thm:main-general}. 
On any disk $|\lambda|<M$ for $M>0$, Lemma~\ref{lem:closed-range} realizes all of them as kernels of one bounded operator on one Banach space $X_\nu$, with closed range throughout
the disk. 
For $\lambda\in \C$, let $n_p(\lambda)=\dim(K_p(\lambda))$.

Apply M\"uller's punctured-neighborhood theorem
\cite[Theorems~1.6 and~1.8]{Muller2000} to the operator complex
\[0\longrightarrow X_\nu \xrightarrow{\,T_p-\lambda\,} X_\nu. \]
It follows that the family $T_p-\lambda$ is regular in a punctured neighborhood of every point, and that its nullity there is
constant and no larger than its value at the center.
Consequently, there are an integer $n_0\ge0$ and a set $\Delta\subset\C$ of \emph{exceptional points}, with no finite accumulation point, such that
\begin{equation}\label{e:generic-nullity}
n_p(\lambda)=n_0\quad \text{ when } \lambda\not\in \Delta \qquad\text{ and } \qquad 
n_p(\lambda)>n_0 \quad\text{ when } \lambda\in \Delta. 
\end{equation}
Indeed, the constant values of $n_p$ on two intersecting punctured
neighborhoods agree, since the intersection contains a point other than
either center. 
Since $\lambda\in \C$, which is connected, 
the constant is global.
Any $\lambda_0\in \C$ has a neighborhood containing no exceptional point other than possibly $\lambda_0$ itself. 
Although the Banach space used in the argument depends on $M$, \eqref{e:full-kernel} identifies its kernel with $K_p(\lambda)$ whenever $|\lambda| < M$. 
Thus the different choices of $M$ give the same value of $n_p(\lambda)$ wherever both apply. 

Fix $\lambda_0\in \C$ such that $n_p(\lambda)$ is constant in a neighborhood of $\lambda_0$. 
We show that the family $T_p-\lambda$, by varying $\lambda$, is regular there. 
To see this directly, write $X_\nu=K_p(\lambda_0)\oplus Y$.
By Lemma \ref{lem:closed-range}, the range of $T_p-\lambda$ is closed and it makes $(T_p-\lambda_0)|_Y$ bounded below, and the same lower bound, halved, holds nearby. 
For $x=u+v\in K_p(\lambda)$ with $u\in K_p(\lambda_0)$ and $v\in Y$, applying $T_p-\lambda$ gives $(T_p-\lambda)v=(\lambda-\lambda_0)u$, so the lower bound on $Y$ gives $\|v\|\ll|\lambda-\lambda_0|\|u\|$. 
The projection onto $K_p(\lambda_0)$
is injective on $K_p(\lambda)$ and, by equal dimensions, surjective.
Thus, as $\lambda\to \lambda_0$, every vector of norm at most one in either $K_p(\lambda)$ or $K_p(\lambda_0)$ can be approximated by a vector in the other space, with the error tending uniformly to zero. 
The regularity of the family $T_p-\lambda$ follows from \cite[Theorem~1.1]{Muller2000}.

We now prove that $n_0=0$. 
Suppose for contradiction that $n_0>0$. Since $T_p$ and $T_q$ commute, $K_p(\lambda)$ is stable under $T_q$. For $\lambda\notin\Delta$ define
\[ d(\lambda)=\det(T_q|K_p(\lambda)).\]
We verify that $d$ extends to an entire function, using only local bases away from isolated exceptional points.
Fix a center $\lambda_0$, choose $M>|\lambda_0|+1$, and let $\nu$ be as in Lemma~\ref{lem:closed-range} for this $M$, so that every $T_p-\lambda$ with $|\lambda|<M$ acts on the same space $X_\nu$.
On a sufficiently small punctured disk the family $T_p-\lambda$ is regular. 
At any point $\lambda'$ of this punctured disk,
\cite[Theorem~1.2]{Muller2000} extends a basis of $K_p(\lambda')$ to
holomorphic kernel vectors $x_1(\lambda), \ldots, x_{n_0}(\lambda)$
near $\lambda'$. 
Let $x_1^*,\ldots,x_{n_0}^*\in X_\nu^*$ be the dual basis to $x_1(\lambda'),\ldots,x_{n_0}(\lambda')$, extended to $X_\nu$ by Hahn--Banach.
Then
\[A(\lambda) =(x_j^*(x_i(\lambda)))_{j,i},\qquad
B(\lambda)=(x_j^*(T_qx_i(\lambda)))_{j,i}\]
are holomorphic, $A(\lambda')=I$, and $A^{-1}B$ represents $T_q|K_p(\lambda)$ near $\lambda'$. 
Hence $d=\frac{\det B}{\det A}$ is holomorphic there. 
Also
\[|d(\lambda)| \le \|T_q\|_{X_\nu}^{n_0},\]
because every eigenvalue of the restriction is bounded by the ambient operator norm. 
Thus the singularity at $\lambda_0$ is removable. 
If $\lambda_0\notin\Delta$, this extension has the prescribed determinant value, because $T_p-\lambda$ is regular at $\lambda_0$ by the previous paragraph.
The local extensions therefore agree and give an entire $d$.
No choice of a basis at a jumping kernel is needed.

For $\lambda\notin\Delta$, let $\mu_1,\ldots,\mu_{n_0}$ be the eigenvalues of $T_q$ on the space $K_p(\lambda)$, counted with algebraic multiplicity.
Each has a nonzero eigenvector in $K_p(\lambda)$, so
Lemma~\ref{lem:two-prime} applies to each pair $(\lambda,\mu_i)$.
With $\alpha=\log q/\log p$ it follows that, for large $|\lambda|$,
\begin{equation}\label{e:det-growth}
 c_{p,q}^{n_0}|\lambda|^{n_0\alpha}
 \le |d(\lambda)|\le
 C_{p,q}^{n_0}|\lambda|^{n_0\alpha}.
\end{equation}
Continuity of $d(\lambda)$ extends these bounds to the exceptional points as well.
The upper bound and Cauchy's estimates imply that $d$ is a
polynomial; the lower bound makes it nonzero.  
Comparing at infinity in \eqref{e:det-growth} gives the degree of the polynomial is $n_0\log q/\log p$, which should be an integer. 
Since $p$ and $q$ are distinct primes, this contradicts $n_0>0$.
Thus $n_0=0$, and by \eqref{e:generic-nullity}, $K_p(\lambda)\ne0$ if and only if $\lambda\in\Delta$.
Since $\Delta$ has no finite accumulation point, it meets every bounded subset of $\C$ in only finitely many points.
\end{proof}

Here discreteness concerns the eigenvalues; no assertion about the rest of the spectrum of $T_p$ on a fixed Banach space is required.

\begin{proof}[Proof of Theorem~\ref{thm:main1}]
Fix a prime $p>\sqrt N$ not dividing $N$. For
$(f_n)\in\E_\sigma(N,k)$, let $\alpha,\beta$ be the roots of
$t^2-f_pt+1$, so
\[\sum_{j\ge0}f_{p^j}t^j =\frac1{(1-\alpha t)(1-\beta t)},\qquad \alpha\beta=1.\]
The bound $f_{p^j}=O_f(p^{j\sigma})$ makes the series holomorphic
for $|t|<p^{-\sigma}$. 
Its rational expression has no canceled
poles, so $|\alpha|,|\beta|\le p^\sigma$. 
Thus $\sigma<0$ is impossible, and for $\sigma\ge0$ we have the uniform bound
$|f_p|\le p^\sigma+p^{-\sigma}$.
There are at most $2^{\omega(N)+1}$ sign systems \eqref{e:Euler-Ramanujan-e}, each satisfying \eqref{e:good-signs}.
By Proposition~\ref{prop:point-discrete}, for each $\vecepsilon$, only finitely many $f_p$ occur in this bounded disk. 
For each one, Theorem~\ref{thm:main-general} gives a finite-dimensional space
$K_p(f_p)=\calH_{\vecepsilon}(N,k;\{p\},(f_p))$. 
By Lemmas~\ref{lem:Euler-Ramanujan} and~\ref{lem:Hecke-Euler}, the function $f$ attached to $(f_n)\in\E_\sigma(N,k)$ lies in $K_p(f_p)$ for the sign system $\vecepsilon$ given by \eqref{e:Euler-Ramanujan-e}. 
The functions attached to $\E_\sigma(N,k)$ therefore lie in a finite sum of finite-dimensional spaces. 
They are linearly independent by Lemma~\ref{lem:linindep}, so there are only finitely many of them.
\end{proof}

\newcommand{\etalchar}[1]{$^{#1}$}
\providecommand{\bysame}{\leavevmode\hbox to3em{\hrulefill}\thinspace}
\providecommand{\MR}{\relax\ifhmode\unskip\space\fi MR }
\providecommand{\MRhref}[2]{\href{http://www.ams.org/mathscinet-getitem?mr=#1}{#2}
}
\providecommand{\href}[2]{#2}


\begin{thebibliography}{BBB{\etalchar{+}}18}

\bibitem[AS03]{AS03}
Jean-Paul Allouche and Jeffrey Shallit, \emph{Automatic sequences}, Cambridge
  University Press, Cambridge, 2003, Theory, applications, generalizations.
  \MR{1997038}

\bibitem[BBB{\etalchar{+}}18]{src}
Sandro Bettin, Jonathan~W. Bober, Andrew~R. Booker, Brian Conrey, Min Lee,
  Giuseppe Molteni, Thomas Oliver, David~J. Platt, and Raphael~S. Steiner,
  \emph{A conjectural extension of {H}ecke's converse theorem}, Ramanujan J.
  \textbf{47} (2018), no.~3, 659--684. \MR{3874812}

\bibitem[BCOZ20]{BCOZ}
Benjamin Bedert, George Cooper, Thomas Oliver, and Pengcheng Zhang,
  \emph{Twisting moduli for {$GL(2)$}}, J. Number Theory \textbf{217} (2020),
  142--162. \MR{4140623}

\bibitem[BKK21]{BKK}
Jakub Byszewski, Jakub Konieczny, and El\.zbieta Krawczyk, \emph{Substitutive
  systems and a finitary version of {C}obham's theorem}, Combinatorica
  \textbf{41} (2021), no.~6, 765--801. \MR{4357430}

\bibitem[BLZ26]{github}
Andrew~R. Booker, Min Lee, and Nina Zubrilina, \emph{Supplementary material for
  ``{T}owards a twistless converse theorem''}, GitHub repository,
  \url{https://github.com/arbooker/twistless_converse}, 2026.

\bibitem[BS26]{BS}
Andrew~R. Booker and Andrew~V. Sutherland, \emph{{$L$}-functions from nothing},
  in preparation, 2026.

\bibitem[CF95]{CF}
J.~B. Conrey and D.~W. Farmer, \emph{An extension of {H}ecke's converse
  theorem}, Internat. Math. Res. Notices (1995), no.~9, 445--463. \MR{1360623}

\bibitem[CFOS07]{CFOS}
J.~B. Conrey, David~W. Farmer, B.~E. Odgers, and N.~C. Snaith, \emph{A converse
  theorem for {$\Gamma_0(13)$}}, J. Number Theory \textbf{122} (2007), no.~2,
  314--323. \MR{2292257}

\bibitem[Cob69]{Cobham}
Alan Cobham, \emph{On the base-dependence of sets of numbers recognizable by
  finite automata}, Math. Systems Theory \textbf{3} (1969), 186--192.
  \MR{250789}

\bibitem[Dim]{Dimitrov}
Vesselin Dimitrov, \emph{Arithmetic criteria for the modularity of a formal
  {D}irichlet series},
  \href{https://www.youtube.com/watch?v=lm9cUtmXt8M}{youtube.com/watch?v=lm9cUtmXt8M},
  33\`emes Journ\'ees Arithm\'etiques.

\bibitem[FKL08]{FKL}
David~W. Farmer, Sally Koutsoliotas, and Stefan Lemurell, \emph{Modular forms
  and {$L$}-functions with a partial {E}uler product}, J. Ramanujan Math. Soc.
  \textbf{23} (2008), no.~2, 105--121. \MR{2432791}

\bibitem[Fur67]{Furstenberg}
Harry Furstenberg, \emph{Disjointness in ergodic theory, minimal sets, and a
  problem in {D}iophantine approximation}, Math. Systems Theory \textbf{1}
  (1967), 1--49. \MR{213508}

\bibitem[Hec36]{Hecke}
Erich Hecke, \emph{{\"Uber die Bestimmung Dirichletscher Reihen durch ihre
  Funktionalgleichung}}, Math. Ann. \textbf{112} (1936), no.~1, 664--699.
  \MR{1513069}

\bibitem[Miy89]{Miyake}
Toshitsune Miyake, \emph{Modular forms}, Springer-Verlag, Berlin, 1989.

\bibitem[M{\"u}l00]{Muller2000}
Vladim{\'i}r M{\"u}ller, \emph{On the punctured neighbourhood theorem}, Journal
  of Operator Theory \textbf{43} (2000), no.~1, 83--95.

\bibitem[Per04]{Perelli}
Alberto Perelli, \emph{A survey of the {S}elberg class of {$L$}-functions.
  {II}}, Riv. Mat. Univ. Parma (7) \textbf{3*} (2004), 83--118. \MR{2128842}

\bibitem[Ste19]{S19}
Raphael~S. Steiner, \emph{Near counterexamples to {W}eil's converse theorem},
  J. Th\'eor. Nombres Bordeaux \textbf{31} (2019), no.~2, 313--321.
  \MR{4030909}

\bibitem[Wei67]{Weil}
Andr\'e Weil, \emph{\"{U}ber die {B}estimmung {D}irichletscher {R}eihen durch
  {F}unktionalgleichungen}, Math. Ann. \textbf{168} (1967), 149--156.
  \MR{207658}

\bibitem[Zha26]{Zhang26}
Pengcheng Zhang, \emph{Twisting-free converse theorems for ${\Gamma}_0(32)$ and
  ${\Gamma}_0(36)$}, 2026,
  \href{https://arxiv.org/abs/2611.00850}{arXiv:2611.00850}.

\end{thebibliography}
\end{document}